\RequirePackage{etex} 
\documentclass[10pt,preprint]{elsarticle}

\usepackage[utf8]{inputenc}
\usepackage[T1]{fontenc}
\usepackage{lmodern}
\usepackage[ngerman, english]{babel}
\usepackage{mathtools}
\usepackage{amsmath}
\usepackage{amssymb}
\usepackage[hidelinks]{hyperref}
\usepackage[usenames,dvipsnames,svgnames,table,x11names,rgb,html]{xcolor}
\usepackage{xspace}
\usepackage{bm}
\usepackage{color}
\usepackage{booktabs}
\usepackage{caption}
\usepackage[]{todonotes}
\usepackage{epstopdf}
\usepackage{tikz}
\usepackage{pgfplots}
\usepackage{etex}
\usepackage{listings}
\usepackage{appendix}
\usepackage{algorithm}
\usepackage{algpseudocode}
\renewcommand{\algorithmiccomment}[1]{\bgroup\quad//~#1\egroup}
\usepackage{lipsum}
\usepackage{siunitx}
\usepackage[super]{nth}

\usetikzlibrary{patterns}
\usetikzlibrary{backgrounds}
\usetikzlibrary{calc}
\usetikzlibrary{arrows.meta}
\usetikzlibrary{decorations.pathreplacing}
\usetikzlibrary{decorations.pathmorphing}

\pgfplotsset{compat=1.12}

\usetikzlibrary{external}

\makeatletter
\renewcommand{\todo}[2][]{\tikzexternaldisable\@todo[#1]{#2}\tikzexternalenable}
\newcommand{\td}[1]{\tikzexternaldisable\@todo[inline,caption={},color=yellow!40]{#1}\tikzexternalenable}
\newcommand{\frage}[2]{\tikzexternaldisable\@todo[inline,caption={},color=violet!20]{\textbf{#1}: #2}\tikzexternalenable}
\makeatother

\definecolor{sbase03}{HTML}{002B36}
\definecolor{sbase02}{HTML}{073642}
\definecolor{sbase01}{HTML}{586E75}
\definecolor{sbase00}{HTML}{445258} 
\definecolor{sbase0}{HTML}{839496}
\definecolor{sbase1}{HTML}{788989} 
\definecolor{sbase2}{HTML}{EEE8D5}
\definecolor{sbase3}{HTML}{FDF6E3}
\definecolor{syellow}{HTML}{B58900}
\definecolor{sorange}{HTML}{CB4B16}
\definecolor{sred}{HTML}{DC322F}
\definecolor{smagenta}{HTML}{D33682}
\definecolor{sviolet}{HTML}{6C71C4}
\definecolor{sblue}{HTML}{268BD2}
\definecolor{scyan}{HTML}{1F7972} 
\definecolor{sgreen}{HTML}{859900}

\newcommand{\diff      } [0] {\ensuremath{\,\mathrm{d}}\xspace}

\newcommand{\vect      } [1] {\bm{#1}}

\newcommand{\bvect     } [1] {\overset{\text{\tiny$\leftrightarrow$}}{\vect{#1}}}

\newcommand{\pderivative}[2] {\frac{\partial #1}{\partial #2}}

\newcommand{\transp    } [0] {\mathsf{T}}

\newcommand{\secref    } [1] {Section~\ref{sec:#1}\xspace}

\newcommand{\secsref   } [1] {Sections~\ref{sec:#1}\xspace}

\newcommand{\figref    } [1] {Fig.~\ref{fig:#1}\xspace}

\newcommand{\eqnref    } [1] {\eqref{eqn:#1}\xspace}

\newcommand{\tabref    } [1] {Table~\ref{tab:#1}\xspace}

\newcommand{\appref    } [1] {\ref{app:#1}\xspace}

\newcommand{\alref     } [1] {Algorithm~\ref{alg:#1}\xspace}

\newcommand{\alsref    } [1] {Algorithms~\ref{alg:#1}\xspace}

\graphicspath{{./figures/}}

\makeatletter
\algrenewcommand\ALG@beginalgorithmic{\small}
\makeatother

\usepackage{subcaption}
\usepackage{float}
\usepackage{placeins}

\usepackage[a4paper]{geometry}

\journal{Journal of Computational Physics}

\addto\captionsenglish{}

\biboptions{sort&compress}

\usepackage[normalem]{ulem}

\definecolor{reviewer2}{HTML}{006400} 
\definecolor{reviewer3}{HTML}{00008B} 
\definecolor{reviewerany}{RGB}{217,95,2}
\newcommand*{\addrev}[2]{{{{#2}}}}
\newcommand*{\addrevtwo}[2]{{{{#2}}}}

\newcommand{\trixi}{\texttt{Trixi.jl}\xspace}
\newcommand{\phix}{q_1}
\newcommand{\phiy}{q_2}
\newcommand{\CFL}{\mathtt{CFL}}
\newcommand{\CFLE}{\mathtt{CFL}_{\mathrm{Eu}}}
\newcommand{\CFLHG}{\mathtt{CFL}_{\mathrm{Gr}}}

\numberwithin{equation}{section}
\allowdisplaybreaks 

\usepackage{stmaryrd}

\newcommand{\eg}[0]{{e.g.\@}\xspace}

\usepackage{lineno}
\linenumbers

\begin{document}

\begin{frontmatter}
  \title{A purely hyperbolic discontinuous Galerkin approach for self-gravitating gas dynamics}

  \author[dep1,dep1a]{Michael Schlottke-Lakemper\corref{cor}}
  \ead{mschlott@math.uni-koeln.de}
  \author[dep2]{Andrew R.\ Winters}
  \author[dep3,dep4]{Hendrik Ranocha}
  \author[dep1,dep1a]{Gregor J.\ Gassner}

  \address[dep1]{Department of Mathematics and Computer Science, University of Cologne, Germany}
  \address[dep1a]{Center for Data and Simulation Science, University of Cologne, Germany}
  \address[dep2]{Department of Mathematics; Applied Mathematics, Link\"{o}ping University, Sweden}
  \address[dep3]{King Abdullah University of Science and Technology (KAUST),
  Computer Electrical and Mathematical Science and Engineering Division (CEMSE),
  Thuwal, Saudi Arabia.}
  \address[dep4]{Present address: Applied Mathematics M\"unster,
  University of M\"unster,
  Germany.}
  \cortext[cor]{Corresponding author.}

  \begin{abstract}
One of the challenges when simulating astrophysical flows with self-gravity is to compute the
gravitational forces. In contrast to the hyperbolic hydrodynamic equations, the gravity field is
described by an elliptic Poisson equation. We present a purely hyperbolic approach by reformulating
the elliptic problem into a hyperbolic diffusion problem, which is solved in pseudotime, 
using the same explicit high-order discontinuous Galerkin method we use for the flow solution. The
flow and the gravity solvers operate on a joint hierarchical Cartesian mesh and are two-way coupled
via the source terms. A key benefit of our approach is that it allows the reuse of existing
explicit hyperbolic solvers without modifications, while retaining their advanced features such as
non-conforming and solution-adaptive grids. By updating the gravitational field in each Runge-Kutta
stage of the hydrodynamics solver, high-order convergence is achieved even in coupled multi-physics
simulations. After verifying the expected order of convergence for single-physics and multi-physics
setups, we validate our approach by a simulation of the Jeans gravitational instability.
Furthermore, we demonstrate the full capabilities of our numerical framework by computing a
self-gravitating Sedov blast with shock capturing in the flow solver and adaptive mesh refinement
for the entire coupled system.

  \end{abstract}

  \begin{keyword}
    discontinuous Galerkin spectral element method
    \sep multi-physics simulation
    \sep adaptive mesh refinement
    \sep compressible Euler equations
    \sep hyperbolic self-gravity
  \end{keyword}
\end{frontmatter}

\section{Introduction}
\label{sec:intro}

Numerical simulations of self-gravitating gas dynamics have become an indispensable tool in the
investigation of astrophysical fluid dynamics, as evidenced by the multitude of publicly available
simulation codes \cite{StoneNorman92, fryxell2000flash, Teyssier02, AlmgrenBecknerEtAl10,
BryanNormanEtAl14, hubber2018gandalf, StoneTomidaEtAl20}.  The gravitational effect of matter on
itself and its surroundings plays an important role in many such flow problems, e.g., for
cosmological structure formation \cite{YangRickerEtAl09, ZuHoneMarkevitchEtAl10}, core-collapse
supernovae \cite{CouchGrazianiEtAl13}, or star formation \cite{LatifZaroubiEtAl10,
FederrathKlessen12}. In non-relativistic simulations, self-gravity is modelled by a Poisson
equation for the Newtonian gravitational potential $\phi$,
\begin{equation}\label{eqn:Poisson}
  -\vec{\nabla}^2\phi = -4\pi G \rho,
\end{equation}
where $G$ is the universal gravitational constant and $\rho$ is the mass density. A particular
challenge for simulating self-gravity is that \eqnref{Poisson} is of elliptic type, i.e., the
solution at each point in space depends on the solution at all other points simultaneously.
So, typically, alternative solution methods than those employed for the hyperbolic hydrodynamics equations are required.
Examples are multigrid methods \cite{Ricker08}, multipole expansion \cite{CouchGrazianiEtAl13}, or
tree-based algorithms \cite{BarnesHut86, WuenschWalchEtAl18}. Some commonalities of these methods 
are their computational expense and that they can be difficult to parallelize due to the global
nature of the problem statement.

In 2007, Nishikawa \cite{Nishikawa07} introduced a new strategy for determining the steady-state
solution to the diffusion equation by rewriting it as a first-order hyperbolic system (FOHS) with a
relaxation time. He noted that when all derivatives with respect to time become zero, the FOHS
reduces to an elliptic equation. Thus, by relaxing the FOHS to steady state, it is
possible to recover the solution for Laplace- and Poisson-type equations.

Based on this analysis, we present a new numerical approach for self-gravitating gas dynamics, where
we follow Nishikawa's ansatz and reformulate \eqnref{Poisson} as a first-order hyperbolic system.
That is, we determine the gravitational potential of a given density distribution as the
steady-state solution to a FOHS with the appropriate Poisson source term of
\eqnref{Poisson}. A key benefit of this strategy is that it allows us to use
the \emph{same explicit discontinuous Galerkin scheme} for both gravity and gas dynamics,
yielding a comparatively simple \emph{high-order multi-physics approach} for hydrodynamics with
self-gravity.  Instead of requiring special treatment of the elliptic equation for gravity, it is
sufficient to set up a hyperbolic solver for each physical system, which are coupled via the
respective source terms. As an additional advantage, this approach enables us to exploit advanced
features of existing solvers for hyperbolic equations without further modifications, such as
local mesh refinement and solution-adaptive grids.

In Nishikawa's original paper \cite{Nishikawa07}, he reformulated the diffusion equation,
\begin{equation}\label{eqn:diffusion}
  u_t - \nu \vec{\nabla}^2 u = 0,
\end{equation}
as the first-order hyperbolic system,
\begin{equation}\label{eqn:fohs}
  \begin{aligned}
    u_t - \nu (q_{1,x} + q_{2,y}) &= 0,\\
    q_{1,t} - \frac{u_x}{T_r} &=-\frac{q_1}{T_r},\\
    q_{2,t} - \frac{u_y}{T_r} &=-\frac{q_2}{T_r}.\\
  \end{aligned}
\end{equation}
Here, $\nu$ is the diffusion coefficient, $T_r$ is the relaxation time, and $q_1,q_2$ are auxiliary
variables. The first approach to rewrite \eqnref{diffusion} as a hyperbolic problem dates back
to 1958, when Cattaneo \cite{cattaneo1958form} introduced the hyperbolic heat equations.
Later, Nagy et al.\ \cite{NagyOrtizEtAl94} showed that solutions of \eqnref{fohs} converge uniformly to solutions of
\eqnref{diffusion} as $T_r$ goes to zero. However, the system
\eqnref{fohs} becomes stiff for very small values of $T_r$ and, thus, prohibitively expensive to solve with an explicit time
integration scheme. Implicit schemes, on the other hand, are generally more difficult to
parallelize and can suffer from reduced solution accuracy due to the stiff source term
\cite{LevequeYee90}.

It was Nishikawa who realized that the key property of \eqnref{fohs} is that it is
equivalent to the original equation \eqnref{diffusion} at the steady state for \emph{any} value
of $T_r$ \cite{Nishikawa07}, thereby avoiding the stiffness problem. That is, if all derivatives with
respect to time are zero, we obtain
\begin{equation}\label{eqn:laplace}
  \left.\begin{aligned}
    - \nu (q_{1,x} + q_{2,y}) &= 0\\
    u_x &= q_1\\
    u_y &= q_2\\
  \end{aligned}\,\right\}
  \quad \to \quad - \nu \vec{\nabla}^2 u = 0.
\end{equation}
Thus by solving the hyperbolic system \eqnref{fohs} to steady state, we can obtain the solution to
an elliptic problem at arbitrary precision. This idea has been successfully applied to develop
finite volume-type methods, e.g., \cite{li2018new,Nishikawa14,ahn2020hyperbolic,chamarthi2018high},
as well as discontinuous Galerkin-type methods, e.g.,
\cite{lou2019reconstructed,mazaheri2016efficient}, that approximate the solution of parabolic
partial differential equations.

Other attempts to reformulate the Poisson equation for gravity to reduce the associated
computational complexity have been made. Black and Bodenheimer \cite{BlackBodenheimer75}
recast the elliptic equation as a parabolic equation to solve it iteratively to steady state with an
alternating-direction implicit scheme \cite{PeacemanRachford55}. While this approach has
been implemented successfully \cite{KrebsHillebrandt83, NormanWinkler86, StoneNorman92}, it
still requires a specialized solver for the gravitational potential.  Hirai et al.\
\cite{HiraiNagakuraEtAl16} proposed to rewrite \eqnref{Poisson} as an inhomogeneous wave
equation, motivated by the hyperbolicity of general relativity. Their hyperbolic system recovers
Newtonian gravity in the limit of infinite propagation speed for gravitation, with all the
associated issues of becoming a stiff problem. In practice, they found that the propagation speed
can be taken relatively small without significantly degrading the solution, as long as it exceeds
the characteristic velocity of the coupled flow problem. While their scheme is far more efficient
than a direct Poisson solver, they also reported that the exact value for the propagation speed is
problem-dependent and that the incurred modeling error strongly depends on the boundary conditions
and domain size.

\linelabel{lne:intro_motive_start}
\label{rev1:r2q7}\label{rev1:r3q1part2}\addrev{any}{We note that restricting to the particular application of self-gravity multi-physics problems yields a slightly different
mindset compared to the goal of developing a general Poisson-type solver. For instance, it is the case that one 
has a very good initial solution guess for the gravity system due to its gradual evolution with respect to the fluid density, and most self-gravity problems in practice require adaptive mesh refinement.}
\addrev{any}{It is within the purview of self-gravitating gas dynamics
that we bring together ideas and apply a hyperbolic diffusion methodology to approximate the solution of a Poisson-type problem
with a high-order discontinuous Galerkin method.}

\addrev{any}{A principal goal of this work is to explore this hyperbolic gravity approach, in
particular, if it retains the high-order accuracy of the underlying approximation and if it offers a
viable alternative 
for self-gravitating applications.}
\addrev{any}{Additionally, we} 
\linelabel{lne:intro_motive_end}
demonstrate that it is possible to reuse existing numerical tools for hyperbolic
problems to compute the gravitational potential.
For the flow field, the compressible Euler
equations are discretized in space by a high-order discontinuous Galerkin method and integrated in
time by an explicit Runge-Kutta scheme. The \emph{same numerical solver} is also used to determine
the gravitational field by advancing the corresponding FOHS to steady
state. In a coupled simulation, the gravitational potential corresponding to the
current density distribution is determined before each Runge-Kutta stage of the
hydrodynamics solver. The 
resulting gravitational forces are then used in the source term during the subsequent
Runge-Kutta stage of the flow solver.  Both solvers operate on a shared hierarchical Cartesian mesh
that can be adaptively refined to match dynamically changing resolution requirements.
We verify the accuracy of the numerical methods by showing high-order convergence for the hydrodynamics
and gravity solvers both in single-physics and multi-physics computations. In addition to a standard
explicit Runge-Kutta scheme, we also discuss a Runge-Kutta scheme optimized for reaching the
steady-state solution of the FOHS for the gravitational potential with a given spatial semi-discretization faster. The suitability of our
approach for applied astrophysics problems is validated by performing coupled
flow--gravity simulations of the Jeans gravitational instability and of a Sedov explosion with
self-gravity. All results are obtained with the open-source simulation framework \trixi
\cite{SchlottkeLakemperGassnerEtAl20}, and the corresponding numerical setups are publicly available
to facilitate reproducibility \cite{SchlottkeLakemperWintersEtAl20}.

The manuscript is organized as follows: In the next section, we will introduce the governing
equations for purely hyperbolic self-gravitational gas dynamics, present the used numerical
methods and discuss how the coupling of gas dynamics and gravity is
achieved. Results for
the individual single-physics solvers and for fully coupled flow-gravity simulations are given in
\secref{results}. Finally, in \secref{conclusions} we summarize our findings.
Further details on the algorithmic and implementation details are provided in \appref{algorithms_implementation}.

\section{Mathematical model and numerical methods}
\label{sec:model_methods}

In this section, we introduce the governing equations for self-gravitating gas dynamics with
the compressible Euler equations and show how to reformulate the Poisson equation for the gravity
potential as a hyperbolic diffusion system. This is followed by an outline of the discontinuous
Galerkin spectral element method (DGSEM) used for the spatial discretization.
After, we present the employed time integration methods and discuss optimized schemes
for the hyperbolic gravity system. Finally, we discuss how multi-physics coupling is achieved 
between the two hyperbolic solvers.

\subsection{Governing equations for self-gravitating gas dynamics}
\label{sec:governing_equations}

In the following, we present the equations for compressible fluids under the influence of a
gravitational potential. First, the compressible Euler equations are given in their standard form
with source terms, proportional to the gravity potential, in the momenta and total energy
\cite{Chandrasekhar1961}:
\begin{equation}\label{eqn:eulerPlusGravity}
\pderivative{}{t}\begin{bmatrix}
\rho\\[0.1cm]
\rho v_1\\[0.1cm]
\rho v_2\\[0.1cm]
E\\[0.1cm]
\end{bmatrix}
+
\pderivative{}{x}\begin{bmatrix}
\rho v_1\\[0.1cm]
\rho v_1^2 + p\\[0.1cm]
\rho v_1 v_2\\[0.1cm]
(E+p)v_1\\[0.1cm]
\end{bmatrix}
+
\pderivative{}{y}\begin{bmatrix}
\rho v_2\\[0.1cm]
\rho v_1v_2\\[0.1cm]
\rho v_2^2 + p\\[0.1cm]
(E+p)v_2\\[0.1cm]
\end{bmatrix}
=
\begin{bmatrix}
0\\[0.1cm]
-\rho\phi_x\\[0.1cm]
-\rho\phi_y\\[0.1cm]
-(\vec{v}\cdot\vec{\nabla}\phi)\rho\\[0.1cm]
\end{bmatrix}.
\end{equation}
Here, $\rho$ is the density, $\vec{v} = (v_1,v_2)^\transp$ are the velocities, and $E$ is the total
energy. The pressure $p$ is determined from the ideal gas law
\begin{equation}
p = (\gamma-1)\left(E - \frac{\rho}{2}(v_1^2 + v_2^2)\right),
\end{equation}
with the heat capacity ratio $\gamma$.

Following Nishikawa's work \cite{Nishikawa07}, we convert the Poisson equation for the gravitational
potential \eqnref{Poisson} into the \textit{hyperbolic gravity equations},
\begin{equation}\label{eqn:hypGrav}
\pderivative{}{\addrev{2}{\tau}}\begin{bmatrix}
\phi\\[0.1cm]
\phix\\[0.1cm]
\phiy\\[0.1cm]
\end{bmatrix}
+
\pderivative{}{x}\begin{bmatrix}
-\phix\\[0.1cm]
-\phi/T_r\\[0.1cm]
0\\[0.1cm]
\end{bmatrix}
+
\pderivative{}{y}\begin{bmatrix}
-\phiy\\[0.1cm]
0\\[0.1cm]
-\phi/T_r\\[0.1cm]
\end{bmatrix}
=
\begin{bmatrix}
-4\pi G \rho\\[0.1cm]
-\phix/T_r\\[0.1cm]
-\phiy/T_r\\[0.1cm]
\end{bmatrix},
\end{equation}
where \addrev{2}{$\tau$ is a pseudotime variable,} \linelabel{lne:pseudotime_1}
the auxiliary variables $(\phix,\phiy)^\transp\approx\vec{\nabla}\phi$ and $T_r$ is the relaxation
time. For a general Poisson problem (with viscosity $\nu$), this relaxation time has the form
\begin{equation}
T_r = \frac{L_r^2}{\nu},
\end{equation}
where $L_r$ is a reference length scale that can be freely chosen. For the gravitational potential
equation \eqnref{Poisson} we have the diffusion constant $\nu=1$. 
\linelabel{lne:Lrscaling_start}
\addrev{2}{The work of Nishikawa \cite{Nishikawa14a} 
found that taking the length scale
\begin{equation}\label{eqn:hypDiffLength}
L_r = \frac{1}{2\pi},
\end{equation}
is optimal on unit square domains. Here, optimality is taken to be} in the sense that the convergence of the 
hyperbolic system \eqnref{hypGrav} to steady state is fastest. 
\label{rev1:r2q1}\addrev{2}{For the numerical results presented in \secref{results}, we found that \eqnref{hypDiffLength}
remained optimal on any square domain and, therefore, it is used throughout the present work.
In more general rectangular or irregular domains, Nishikawa and Nakashima demonstrated that the reference length scaling $L_r$
must be adjusted to avoid erratic convergence behavior of the hyperbolic diffusion system
\linelabel{lne:new_ref}\cite{nishikawa2018dimensional}. 
We note, however, that regardless of how one chooses
$L_r$,} 
\linelabel{lne:Lrscaling_end}
the steady-state solution of \eqnref{hypGrav} is, in fact, the desired solution of the
original Poisson problem \eqnref{Poisson} \cite{Nishikawa07, cattaneo1958form, gomez2010hyperbolic,NagyOrtizEtAl94}.

In single-physics simulations of the compressible Euler equations or the hyperbolic gravity
equations alone, the source terms in the governing equations are determined analytically. When
considering coupled Euler-gravity problems, however, the gravity potential of the hyperbolic gravity
system is used to generate the source term information in \eqnref{eulerPlusGravity}, while the
density of the compressible Euler solution is used for the source term in \eqnref{hypGrav}. Thus,
both systems of equations are connected via two-way coupling through their source terms.

\subsection{High-order discontinuous Galerkin method on hierarchical Cartesian meshes}
\label{sec:numerical_methods}
Next, we give an overview of the nodal discontinuous Galerkin spectral element method
(DGSEM) on hierarchical Cartesian meshes. A full derivation can be found in, e.g., \cite{Kopriva09,
HindenlangGassnerEtAl12, Schlottke-LakemperNiemoellerEtAl19}.  We consider the solution of systems of hyperbolic conservation
laws in two spatial dimensions, which take the general form
\begin{equation}\label{eqn:consLawGeneral}
  \vect{u}_t + \vec{\nabla} \cdot \bvect{f}(\vect{u}) = \vect{s}(\vect{u}),
  \qquad\forall \vec{x} \in \Omega,
\end{equation}
on a square domain $\Omega$. Here $\vect{u} \in \mathbb{R}^m$, where $m$ is the number of
equations, is the state vector of conserved variables, $\bvect{f}
= (\vect{f}_1, \vect{f}_2)^\transp$ is the block vector of the nonlinear fluxes, and $\vect{s}$
denotes a---possibly zero-valued---source term.
We subdivide the domain $\Omega$ into $K$ non-overlapping square elements 
\linelabel{lne:mapping_start}
\label{rev1:r2q8}\addrev{2}{
\begin{equation}
E_k = \left[x_{k,1},x_{k,2}\right]\times[y_{k,1},y_{k,2}],\quad k = 1,\ldots,K,
\end{equation}
and we define $h:=\Delta x = \Delta y$ as the length of the respective
Cartesian element. We transform between the reference element $E_0=[-1,1]^2$ and each element, $E_k$, 
from the mappings
\begin{equation}
 X_k(\xi) = x_{k,1} + \frac{\xi +1}{2}h,\quad  Y_k(\eta) = y_{k,1} + \frac{\eta +1}{2}h,\quad k=1,\ldots,K
\label{eqn:mapping}
\end{equation}
with reference coordinates $\vec{\xi} = (\xi,\eta)^\transp$. Under these mappings \eqnref{mapping},
the conservation law in physical coordinates \eqnref{consLawGeneral} becomes a conservation law in reference coordinates
\begin{equation}
\frac{h^2}{4}\vect{u}_t + \frac{h}{2}\vec{\nabla}_{\xi} \cdot \bvect{f} = \frac{h^2}{4}\vect{s}(\vect{u}).
\end{equation}
We simplify the conservation law in reference coordinates to be}
\linelabel{lne:mapping_end}
\begin{equation}\label{eqn:consLawRef}
  J\vect{u}_t + \vec{\nabla}_{\xi} \cdot {\bvect{f}} = J\vect{s},
\end{equation}
where $J = h/2$ is the one-dimensional Jacobian determinant.

The starting point of the DGSEM is the weak form of the conservation law, for which we
multiply \eqnref{consLawRef} by an appropriate test function $\varphi \in
L_2(\Omega)$ and integrate over the reference element. After integration-by-parts, we obtain the
weak form
\begin{equation}\label{eqn:weakFormCont}
  \int\limits_{{E_0}} J \vect{u}_t \varphi \diff \vec{\xi} + \int\limits_{\partial {E_0}} (\bvect{f} \cdot \vec{n}) \varphi
  \diff S
  - \int\limits_{{E_0}} \bvect{f} \cdot \vec{\nabla}_{\xi} \varphi \diff \vec{\xi} =
  \int\limits_{{E_0}} J \vect{s} \varphi \diff \vec{\xi},
\end{equation}
where $\vec{n}$ is the outer unit normal at the boundary $\partial {E_0}$.
We approximate each component of the state vector by polynomials of degree $N$ in each spatial dimension, which we
represent as $\vect{U}$. The polynomials are written in terms of the Lagrange
basis $\psi_{ij}(\,\vec{\xi}\,) = l_i(\xi) l_j(\eta)$, $i,j=0,\ldots,N$, where the interpolation
nodes are the Legendre-Gauss-Lobatto (LGL) points. Lagrange polynomials of
degree $N$ are also used to approximate the fluxes $\bvect{F} = (\vect{F}_1, \vect{F}_2)^\transp$
and source terms $\vect{S}$. Integrals in \eqnref{weakFormCont} are evaluated discretely by
LGL quadrature such that the interpolation and quadrature nodes are collocated.
To resolve the solution discontinuity at element interfaces, we replace the boundary fluxes by
numerical fluxes $\vect{F}_n^* \approx \bvect{F} \cdot \vec{n}$, see, e.g., \cite{Kopriva09}.
In this work, we use the Harten, Lax, van Leer (HLL) flux \cite{harten1983upstream,Toro09} for the
compressible Euler system and
the local Lax-Friedrichs (LLF) flux \cite{Toro09, lou2019reconstructed} for hyperbolic gravity.
We choose the tensor product basis $\varphi \in \{\psi_{ij}\}_{i,j=0}^N$ to be the $(N+1)^2$
test functions. The final weak form of \eqnref{consLawGeneral} in the standard DGSEM formulation
then reads
\begin{equation}\label{eqn:weakFormDisc}
  \int\limits_{{E_0},N} J \vect{U}_t \psi_{ij} \diff \vec{\xi} + \int\limits_{\partial {E_0},N} \vect{F}_n^* \psi_{ij}
  \diff S
  - \int\limits_{{E_0},N} \bvect{F} \cdot \vec{\nabla}_{\xi} \psi_{ij} \diff \vec{\xi} =
  \int\limits_{{E_0},N} J \vect{S} \psi_{ij} \diff \vec{\xi}, \qquad \forall i,j = 0,\ldots,N.
\end{equation}
Inserting the definitions for the approximate solution, fluxes, and source terms, we obtain the
semi-discrete DG operator, which is integrated in
time with an explicit Runge-Kutta scheme. The stable time step is calculated by
\begin{equation}\label{eqn:max_dt}
  \Delta t = \frac{\CFL}{N + 1} \frac{h}{\Lambda},
\end{equation}
where $\Lambda_\mathrm{Eu} = \max(\mathopen{|}\vec{v}\mathclose{|} + c)$
with the speed of sound $c = \sqrt{\gamma p/\rho}$ for the compressible Euler system and
$\Lambda_\mathrm{Gr} = \sqrt{\nu/T_r}$ for the gravity system.
Suitable time integration schemes will be discussed in the next section.

An alternative split-form DGSEM approximation can be obtained by making use of the
summation-by-parts (SBP) property inherent to the nodal DG scheme on LGL nodes \cite{Gassner13a}.
Applying summation-by-parts once more to \eqnref{weakFormDisc} produces the strong form of the
conservation law \eqnref{consLawRef}. Following the work of Fisher et al.\ and Carpenter et al.\
\cite{FisherCarpenterEtAl13, CarpenterFisherEtAl14}, we introduce a numerical volume flux for the
flux derivative. This yields a split-form DG approximation \cite{GassnerWintersEtAl16,
GassnerWintersEtAl18}, which allows to use symmetric two-point flux functions with additional
desirable properties such as entropy conservation or kinetic energy preservation
\cite{tadmor2003entropy, ismail2009affordable, chandrashekar2013kinetic, ranocha2018thesis}. Here,
we select the numerical flux of Chandrashekar
\cite{Chandrashekar2012}.

The split-form variant is also the basis for the high-order shock capturing approach by Hennemann et
al.\ \cite{HENNEMANN2021109935}, which we utilize in the DG solver for simulating
compressible Euler problems with strong discontinuities. In this approach, each DG element is
divided into $(N+1)^2$ subcells, on which a first-order finite volume (FV) method is used to obtain
the semi-discrete operator.  The final operator is then obtained by blending the DG operator in each
element with the FV operator based on the energy content in the highest modes, while retaining the
discrete entropic property.

\begin{figure}[!htp]
  \centering
  \begin{subfigure}[t]{0.49\textwidth}
    \centering
    \includegraphics[width=57mm]{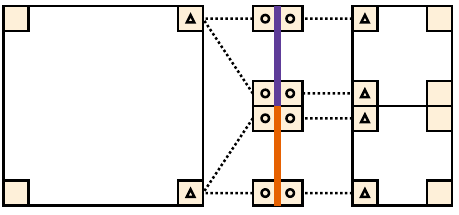}
    \caption{Mortar method: non-conforming nodes
    on element surfaces (square with triangle) become conforming nodes on mortar
    elements (square with circle).}
    \label{fig:mortar}
  \end{subfigure}
  \hfill
  \begin{subfigure}[t]{0.49\textwidth}
    \centering
    \includegraphics[width=69mm]{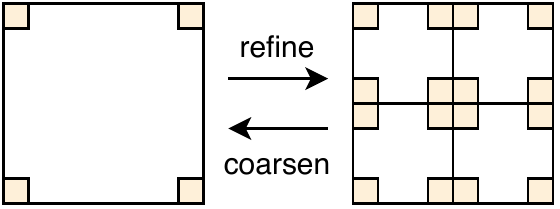}
    \caption{Adaptive mesh refinement: solution transfer between coarse and refined cells.}
    \label{fig:amr}
  \end{subfigure}
  \caption{Illustration of the mortar element method and transfer operators for adaptive mesh
  refinement.}
  \label{fig:mortar_amr}
\end{figure}

To support non-conforming elements created by local mesh refinement, we employ the mortar method
\cite{Kopriva96, KoprivaWoodruffEtAl02}. In this approach, mortar surfaces are inserted at
interfaces with a coarse element adjacent to multiple refined elements, see \figref{mortar}. The
solution values at non-conforming element interfaces are interpolated to the mortar, where the
surface flux is calculated at conforming node locations. Then, the flux values are discretely
projected back to the non-conforming element interfaces. A similar approach is used for adaptive mesh
refinement (AMR), see \figref{amr}. During refinement, the solution on the coarse element is
interpolated to the LGL nodes on the four refined elements. Conversely, during coarsening the
solution on the four refined elements is projected onto the coarse element. These AMR procedures in
two dimensions correspond exactly to the algorithms used for a mortar element approach in three
dimensions \cite{bui2012analysis}. Therefore, both the mortar method and the AMR technique are fully conservative. A detailed
derivation of the DGSEM on non-conforming hierarchical Cartesian meshes can be found in, e.g.,
\cite{Schlottke-LakemperNiemoellerEtAl19}.

\subsection{Time integration schemes}\label{sec:3STime}

For the compressible Euler solver, we use the fourth-order, five-stage low-storage scheme CK45
by Carpenter and Kennedy \cite{CarpenterKennedy94} with CFL-based step size control \eqnref{max_dt}.
This and other common time integration methods can also be used to advance the
hyperbolic diffusion part to steady state.
To improve the performance of integration in pseudotime \cite{loppi2019locally}
and reduce the sensitivity to user-chosen parameters \cite{ranocha2021optimized},
locally adaptive error-based time step control
could be employed and the Runge-Kutta schemes can
be optimized for the given spatial semi-discretizations to be able to take
bigger time steps \cite{vermeire2020optimal}.

Here, we briefly demonstrate the last approach. 
\linelabel{lne:timestep_opt_start}
\label{rev1:r2q3}\addrev{2}{The motivation behind this strategy is
to grant the scheme
an ability to take larger explicit time steps, which leads to faster steady state
convergence for the numerical examples considered in \secref{results}.}
\linelabel{lne:timestep_opt_end}
Following \cite{parsani2013optimized}, we computed the spectrum of the standard DG operator
using a local Lax-Friedrichs interface flux and polynomials of degree $N = 3$ for
the hyperbolic diffusion problem \eqref{eqn:fohs}. The convex hull of this spectrum was used as
input to optimize the stability polynomial of an explicit, second-order accurate five-stage
Runge-Kutta method using the algorithm of \cite{ketcheson2013optimal}. Finally,
a low-storage scheme of the 3S* class \cite{ketcheson2010runge} was constructed
by minimizing
the principal truncation error, given the optimized stability polynomial as
constraint. For this, we used the algorithms implemented in RK-Opt \cite{RK-Opt},
which are based on MATLAB \cite{matlabR2019b}. We used NodePy \cite{nodepy}
to verify the desired properties of the resulting five-stage second-order method RK3S*
with low-storage coefficients listed in \tabref{3Sstar15opt}. The minimum-storage
implementation \cite{ketcheson2010runge} is shown in \alref{3Sstar}.
Applying RK3S* instead of CK45 increases the performance and does not influence
the accuracy negatively, as shown in \secsref{coupleConv} and \ref{sec:jeans}
below. We also optimized first- and third-order accurate schemes analogously
but these did not improve performance as much, e.g., for the test described in
\secref{jeans}.

\begin{table}
\caption{Minimum-storage coefficients \cite{ketcheson2010runge} of the second-order
         accurate explicit low-storage Runge-Kutta method RK3S* optimized for
         hyperbolic diffusion.}
\label{tab:3Sstar15opt}
\centering
\begin{small}%
\begin{tabular}{@{}cccc@{}}
  \toprule
  $i$ & $\gamma_{1,i}$ & $\gamma_{2,i}$ & $\gamma_{3,i}$ \\
  \midrule
  1 & \verb| 0.0000000000000000E+00| & \verb| 1.0000000000000000E+00| & \verb| 0.0000000000000000E+00| \\
  2 & \verb| 5.2656474556752575E-01| & \verb| 4.1892580153419307E-01| & \verb| 0.0000000000000000E+00| \\
  3 & \verb| 1.0385212774098265E+00| & \verb|-2.7595818152587825E-02| & \verb| 0.0000000000000000E+00| \\
  4 & \verb| 3.6859755007388034E-01| & \verb| 9.1271323651988631E-02| & \verb| 4.1301005663300466E-01| \\
  5 & \verb|-6.3350615190506088E-01| & \verb| 6.8495995159465062E-01| & \verb|-5.4537881202277507E-03| \\
  \toprule
  $i$ & $\delta_i$ & $\beta_{i}$ & $c_i$ \\
  \midrule
  1 & \verb| 1.0000000000000000E+00| & \verb| 4.5158640252832094E-01| & \verb| 0.0000000000000000E+00| \\
  2 & \verb| 1.3011720142005145E-01| & \verb| 7.5974836561844006E-01| & \verb| 4.5158640252832094E-01| \\
  3 & \verb| 2.6579275844515687E-01| & \verb| 3.7561630338850771E-01| & \verb| 1.0221535725056414E+00| \\
  4 & \verb| 9.9687218193685878E-01| & \verb| 2.9356700007428856E-02| & \verb| 1.4280257701954349E+00| \\
  5 & \verb| 0.0000000000000000E+00| & \verb| 2.5205285143494666E-01| & \verb| 7.1581334196229851E-01| \\
  \bottomrule
\end{tabular}
\end{small}%
\end{table}

\begin{algorithm}
  \begin{algorithmic}
    \State $S_1 \gets u^n, S_2 \gets 0, S_3 \gets u^n$
    \ForAll{$i \in \{1, \dots, s\}$}
      \State $S_2 \gets S_2 + \delta_i S_1$
      \State $S_1 \gets \gamma_{1,i} S_1 + \gamma_{2,i} S_2 + \gamma_{3,i} S_3 + \beta_{i} \Delta t f(t_n + c_i \Delta t, S_1)$
    \EndFor
    \State $u^{n+1} \gets S_1$
  \end{algorithmic}
  \caption{Minimum-storage implementation of one step of a 3S* method with $s$
           stages applied to the ODE $u_t(t) = f(t, u(t))$.}
  \label{alg:3Sstar}
\end{algorithm}

\subsection{Multi-physics coupling}\label{sec:multiphysics_coupling}

The results presented in \secref{results} were obtained with
\trixi \cite{SchlottkeLakemperGassnerEtAl20},
an open-source numerical simulation framework\footnote{\trixi:
\href{https://github.com/trixi-framework/Trixi.jl}{https://github.com/trixi-framework/Trixi.jl}} for
hyperbolic conservation laws developed by the authors. \trixi is based on a
modular architecture, where all components are only loosely coupled with each other such that it is
easy to extend or replace existing functionality.

\linelabel{lne:multiphys_start}
\addrev{3}{A complete description of the control flow and algorithmic structure
of \trixi for single-physics simulation setups on hierarchical Cartesian meshes is provided in \appref{algorithms_implementation}.}
\addrev{3}{The extension to a coupled flow-gravity simulation is straightforward, since
the compressible Euler hydrodynamic solver and the hyperbolic gravity solver both use a high-order DG spatial approximation.
Therefore, one can instantiate two instances of the DG solver that operate on the same quadtree mesh. We note that only the solution 
information from the compressible Euler equations is used to control the quadtree mesh 
and its possible adaptation via AMR (see \appref{algorithms_implementation} for details). The gravity solver, however, is passively adapted to match the new
quadtree mesh such that both solvers continue to use the same spatial discretization. We highlight the ease with which one
can adapt to multi-physics simulations is due to the treatment of the gravity equations in a hyperbolic fashion and the ability to directly reuse tools from
an existing single-physics DG architecture to approximate the solution of the gravity potential and its derivatives \eqnref{hypGrav}.}

\begin{figure}[htbp]
  \centering
  \includegraphics[width=\textwidth]{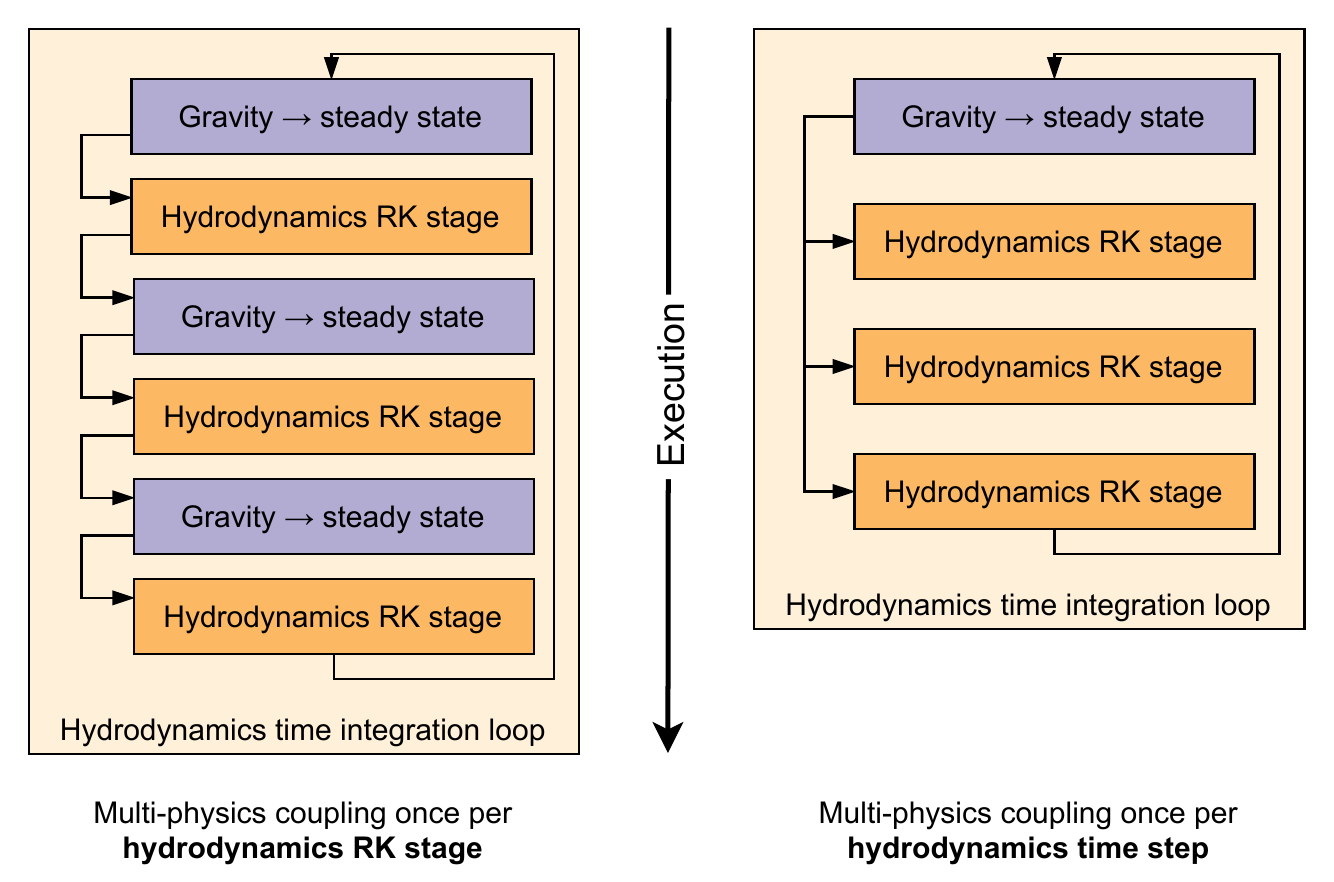}
  \caption{(left) Flowchart for a flow-gravity simulation with coupling once per Runge-Kutta (RK)
  \emph{stage} of the hydrodynamics solver (right) Flowchart for a flow-gravity simulation with
  coupling once per \emph{time step} of the hydrodynamics solver.}
  \label{fig:trixi_flowchart}
\end{figure}

\addrev{3}{The key component of the multi-physics simulation is the exchange of information between the hydrodynamics and hyperbolic gravity solvers via their source terms. 
This coupling introduces a crucial algorithmic choice within the \textit{time integration loop} of the hydrodynamics solver. 
One can choose to evolve the hyperbolic gravity variables to steady state in pseudotime
either within every \emph{stage} of the Runge-Kutta scheme of the hydrodynamics solver (left part of \figref{trixi_flowchart}) or once every
Runge-Kutta time \emph{step} of the hydrodynamics solver (right part of
\figref{trixi_flowchart}).
This choice impacts the overall solution accuracy as we demonstrate with results given in
\secsref{coupleConv} and \ref{sec:jeans}.}
\linelabel{lne:multiphys_end}

\linelabel{lne:new_algs_start}
\addrev{2}{In \alsref{gravEveryStage} and \ref{alg:gravEveryStep} we provide pseudocode for these two coupling strategies. With either
coupling variant, 
we insert a pseudotime integration loop to evolve the hyperbolic gravity solution variables to a new
steady state into the time integration loop of the hydrodynamics solver. As discussed in \secref{3STime},
the time integration scheme for the gravity system can differ from the time integration scheme of
the hydrodynamics system. 
For the integration in pseudotime, the previous gravity solution is used as the initial guess,} 
\linelabel{lne:new_algs_end}
and the source terms as given in \eqnref{hypGrav} are
determined from the density field of the hydrodynamics solver. After the gravity field is converged to steady state
according to the specified residual threshold, the compressible Euler solver is executed accordingly, with the source terms as
given in \eqnref{eulerPlusGravity} computed from the solution of the gravity solver. As a final note, since both solvers operate on 
the same mesh with the same polynomial order, no interpolation is required between the hydrodynamics and gravity solutions.

\begin{algorithm}
  \begin{algorithmic}
    \While{$t < t_{\text{final}}$  }
       \State $\Delta t_{\text{Euler}} \gets$ Value from \eqnref{max_dt} for hydrodynamics solver
       \ForAll{$i \in \{1, \dots, s\}$}
          \State $t_{\text{stage}} \gets t + c_i\Delta t_{\text{Euler}}$
           \While{$\phi_\tau > \mathtt{tol}$}\Comment{Gravity loop in pseudotime $\tau$}
             \State $\Delta t_{\text{Gravity}} \gets$ Value from \eqnref{max_dt} for hyperbolic gravity solver
             \State Update $\vect{u}_{\text{Gravity}}$ using density $\rho$ at $t_{\text{stage}}$ in \eqnref{hypGrav} with gravity solver RK method
          \EndWhile
          \State Update $\vect{u}_{\text{Euler}}$ using $q_1, q_2$ at $t_{\text{stage}}$ in \eqnref{eulerPlusGravity} with hydrodynamics solver RK method
       \EndFor
       \State $t \gets t + \Delta t_{\text{Euler}}$
    \EndWhile
  \end{algorithmic}
  \caption{\addrev{2}{Procedure to couple and update hyperbolic gravity variables once every Runge-Kutta \textbf{stage}
  in the compressible Euler time integration loop.}}
  \label{alg:gravEveryStage}
\end{algorithm}

\begin{algorithm}
  \begin{algorithmic}
    \While{$t < t_{\text{final}}$  }
        \While{$\phi_\tau > \mathtt{tol}$}\Comment{Gravity loop in pseudotime $\tau$}
          \State $\Delta t_{\text{Gravity}} \gets$ Value from \eqnref{max_dt} for hyperbolic gravity solver
          \State Update $\vect{u}_{\text{Gravity}}$ using density $\rho$ at current time $t$ in \eqnref{hypGrav} with gravity solver RK method
       \EndWhile
       \State $\Delta t_{\text{Euler}} \gets$ Value from \eqnref{max_dt} for hydrodynamics solver
       \ForAll{$i \in \{1, \dots, s\}$}
          \State $t_{\text{stage}} \gets t + c_i\Delta t_{\text{Euler}}$
          \State Update $\vect{u}_{\text{Euler}}$ using $q_1, q_2$ at fixed $t$ in \eqnref{eulerPlusGravity} with hydrodynamics solver RK method
       \EndFor
       \State $t \gets t + \Delta t_{\text{Euler}}$
    \EndWhile
  \end{algorithmic}
  \caption{\addrev{2}{Procedure to couple and update hyperbolic gravity variables once every Runge-Kutta \textbf{time step} 
  in the compressible Euler time integration loop.}}
  \label{alg:gravEveryStep}
\end{algorithm}

\section{Numerical results}
\label{sec:results}

In this section, we present numerical tests to verify and validate the multi-physics implementation
of the DG solver in \trixi \cite{SchlottkeLakemperGassnerEtAl20}. We begin with a demonstration, in
\secref{num_verif}, of the high-order accuracy for both single- and multi-physics problems. Next, in
\secref{num_appl}, we simulate two example problems for compressible, self-gravitating flows. All
parameter files required for reproducing the results are also available online
\cite{SchlottkeLakemperWintersEtAl20}.

The time integration of the compressible Euler portion of the multi-physics uses the explicit five stage, fourth order low-storage Runge-Kutta (CK45) scheme of Carpenter and Kennedy~\cite{CarpenterKennedy94}, where a stable time step is computed according to the adjustable coefficient ${\CFLE}\in(0,1]$, e.g., \cite{gassner2011}. To run the hyperbolic diffusion equation system to steady state \addrev{2}{in pseudotime}, \linelabel{lne:pseudotime_2}
we employ either CK45 or RK3S* from \secref{3STime}. A stable time step for the hyperbolic gravity problem is computed with a separate adjustable ${\CFLHG}$ coefficient \cite{lou2019reconstructed,Nishikawa14}. Further, a prescribed tolerance (\texttt{tol}) is set for a given problem to determine when steady state for the hyperbolic gravity solver is reached numerically. A solution is deemed converged when the magnitude of the semi-discrete DG operator for \addrev{2}{$\phi_{\tau}$} \linelabel{lne:pseudotime_3}
is reduced below the
prescribed tolerance in the discrete $L^{\infty}$ norm at the LGL nodes.

\subsection{Verification of single and multi-physics solvers}\label{sec:num_verif}

First, we show that the single-physics implementations for the compressible Euler equations and the hyperbolic diffusion equations are high-order accurate in \secsref{euler_eoc} and \ref{sec:hypdiff_eoc}, respectively.
Then, we demonstrate the high-order accuracy of the coupled, multi-physics solver in \secref{coupleConv}.
For these investigations, we use the standard DGSEM as described in \secref{numerical_methods}.
The explicit time step is selected by setting ${\CFLE}={\CFLHG}=0.5$ such that spatial errors in the
approximation are dominant. The tolerance to
define steady state for the hyperbolic gravity equations is taken as $\mathtt{tol}=10^{-10}$.

\subsubsection{Compressible Euler solver}\label{sec:euler_eoc}

We verify the high-order spatial accuracy for the DG approximation of the compressible Euler equations with the method of manufactured solutions. To do so, consider the system \eqnref{eulerPlusGravity} governing ideal gas dynamics \textit{without} gravitational source terms, i.e., $\phi(x,y)\equiv 0$.

The domain is $\Omega=[0,2]^2$ with periodic boundary conditions and $\gamma=2$. The solution for this test case has the form
\begin{equation}\label{eqn:EulerManufac}
\rho = 2 + \frac{1}{10}\sin(\pi(x+y-t)),\quad v_1 = v_2 = 1,\quad p = \frac{1}{\pi}\rho^2.
\end{equation}
A significant advantage of this ansatz is that it is symmetric and spatial derivatives cancel with temporal derivatives as
\begin{equation}
\rho_x = \rho_y = -\rho_t.
\end{equation}
The manufactured solution produces an additional residual term that reads
\begin{equation}
 \vect{u}_t + \vec{\nabla} \cdot \bvect{f}(\vect{u}) =
\begin{bmatrix}
 \frac{1}{10} \pi \cos(\pi(x+y-t))\\[0.1cm]
 \frac{1}{10} \pi \cos(\pi(x+y-t))\left[1 + \frac{2}{\pi}\left(2 + \frac{1}{10}\sin(\pi(x+y-t))\right)\right]\\[0.1cm]
 \frac{1}{10} \pi \cos(\pi(x+y-t))\left[1 + \frac{2}{\pi}\left(2 + \frac{1}{10}\sin(\pi(x+y-t))\right)\right]\\[0.1cm]
\frac{1}{10} \pi \cos(\pi(x+y-t))\left[1 + \frac{6}{\pi}\left(2 + \frac{1}{10}\sin(\pi(x+y-t))\right)\right]\\[0.1cm]
\end{bmatrix}.
\end{equation}
We run the manufactured solution test case with $T=1.0$ as the final time. For the computations we use uniform Cartesian meshes with a varying number of elements.

To investigate the accuracy of the DG approximation, we use two polynomial orders $N=3$ and $N=4$. We compute the discrete $L^2$ errors in the conservative variables with LGL quadrature over the entire domain at the final time for different mesh resolutions. For the time integration, we select CK45. We present the experimental order of convergence (EOC) in \tabref{EulerEOC} for increasing mesh resolution and the two selected polynomial orders. The results confirm the expected theoretical order of convergence $N+1$ for the DG method, e.g., \cite{Kopriva09}.
\begin{table}[!htb]
    \caption{Convergence test for the compressible Euler equations with manufactured solution \eqnref{EulerManufac} for two polynomial orders.}
    \label{tab:EulerEOC}
    \begin{subtable}{0.5\linewidth}
      \centering
      \begin{scriptsize}
        \caption{$N=3$}
         \begin{tabular}{@{}lcccc@{}}
         \toprule
         $K$ & $L^2(\rho)$ & $L^2(\rho v_1)$ & $L^2(\rho v_2)$ & $L^2(E)$ \\
         \midrule
           $4^2$ & 1.74E-04 & 3.37E-04 & 3.37E-04 & 6.10E-04 \\
           $8^2$ & 1.72E-05 & 2.33E-05 & 2.33E-05 & 4.38E-05 \\
         $16^2$ & 9.64E-07 & 1.39E-06 & 1.39E-06 & 2.62E-06 \\
         $32^2$ & 6.31E-08 & 8.80E-08 & 8.80E-08 & 1.65E-07 \\\midrule
         avg. EOC & 3.81 & 3.97 & 3.97 & 3.95 \\\bottomrule
        \end{tabular}
        \end{scriptsize}
    \end{subtable}%
    \quad
    \begin{subtable}{0.5\linewidth}
      \centering
      \begin{scriptsize}
        \caption{$N=4$}
         \begin{tabular}{@{}lcccc@{}}
         \toprule
         $K$ & $L^2(\rho)$ & $L^2(\rho v_1)$ & $L^2(\rho v_2)$ & $L^2(E)$ \\
         \midrule
           $4^2$ & 1.72E-05 & 2.68E-05 & 2.68E-05 & 4.95E-05 \\
           $8^2$ & 6.82E-07 & 8.92E-07 & 8.92E-07 & 1.68E-06 \\
         $16^2$ & 1.86E-08 & 2.58E-08 & 2.58E-08 & 4.69E-08 \\
         $32^2$ & 6.14E-10 & 8.18E-10 & 8.18E-10 & 1.48E-09 \\\midrule
         avg. EOC & 4.92 & 5.00 & 5.00 & 5.01 \\\bottomrule
        \end{tabular}
        \end{scriptsize}
    \end{subtable}
\end{table}

\subsubsection{Hyperbolic diffusion solver}\label{sec:hypdiff_eoc}

Next, we verify the accuracy of the DG implementation of the first-order hyperbolic diffusion system used to approximate the solution of Poisson's equation. For this we consider the general Poisson problem
\begin{equation}\label{eqn:truePoisson}
-\vec{\nabla}^2\phi = f(x,y),
\end{equation}
which can be converted into a hyperbolic system \cite{Nishikawa07,Nishikawa14a,lou2019reconstructed} analogous to \eqnref{hypGrav}
\begin{equation}\label{eqn:generalFOHS}
\pderivative{}{\addrev{2}{\tau}}\begin{bmatrix}
\phi\\[0.1cm]
\phix\\[0.1cm]
\phiy\\[0.1cm]
\end{bmatrix}
+
\pderivative{}{x}\begin{bmatrix}
-\phix\\[0.1cm]
-\phi/T_r\\[0.1cm]
0\\[0.1cm]
\end{bmatrix}
+
\pderivative{}{y}\begin{bmatrix}
-\phiy\\[0.1cm]
0\\[0.1cm]
-\phi/T_r\\[0.1cm]
\end{bmatrix}
=
\begin{bmatrix}
f(x,y)\\[0.1cm]
-\phix/T_r\\[0.1cm]
-\phiy/T_r\\[0.1cm]
\end{bmatrix}.
\end{equation}
We consider the domain $\Omega=[0,1]^2$ and take the solution and forcing function in \eqnref{truePoisson} to be
\begin{equation}\label{eqn:PoissonManufac}
\phi(x,y) = 2 + 2\cos(\pi x)\sin(2\pi y)\quad\text{ and }\quad f(x,y) = 10\pi^2\cos(\pi x)\sin(2\pi y).
\end{equation}
Analytical expressions for the auxiliary variables $\phix$ and $\phiy$ are then determined by differentiation
\begin{equation}
\phix(x,y) = -2\pi\sin(\pi x)\sin(2\pi y)\quad\text{ and }\quad \phiy(x,y) = 4\pi\cos(\pi x)\cos(2\pi y).
\end{equation}
The boundary conditions are Dirichlet in the $x$-direction and periodic in the $y$-direction.
In the limit of steady state, the hyperbolic system \eqnref{generalFOHS} recovers the solution for
the Poisson equation \eqnref{truePoisson} \cite{Nishikawa07,cattaneo1958form}.

We choose two polynomial orders $N=3$ and $N=4$ to demonstrate the high-order accuracy of the DG solver applied to the hyperbolic diffusion system, and use the RK3S* scheme for time integration. The discrete $L^2$ errors for the variables $\phi$, $\phix$ and $\phiy$ are computed on uniform Cartesian meshes of increasing resolution. The EOCs given in \tabref{HypDiffEOC} are the expected optimal order of $N+1$ in all variables \cite{lou2019reconstructed}. It is interesting to note that the DG approximation for \eqnref{generalFOHS} provides a high-order accurate approximation to $\phi$ as well as its gradient. This is convenient because an accurate approximation of these gradient values is needed for the gravitational coupling terms in \eqnref{eulerPlusGravity}.
\begin{table}[!htb]
    \caption{Convergence test for the hyperbolic diffusion form of the Poisson equation with manufactured solution \eqnref{PoissonManufac} for two polynomial orders. The approximation is high order for the variable $\phi$ as well as its gradient.}
    \label{tab:HypDiffEOC}
    \begin{subtable}{0.5\linewidth}
      \centering
        \begin{scriptsize}
        \caption{$N=3$}
         \begin{tabular}{@{}lccc@{}}
         \toprule
         $K$ & $L^2(\phi)$ & $L^2(\phix)$ & $L^2(\phiy)$\\
         \midrule
          $4^2$ & 3.15E-03 & 1.24E-02 & 2.19E-02\\
          $8^2$ & 2.26E-04 & 8.83E-04 & 1.50E-03\\
         $16^2$ & 1.50E-05 & 5.51E-05 & 9.68E-05\\
         $32^2$ & 9.65E-07 & 3.32E-06 & 6.14E-06\\\midrule
         avg. EOC & 3.89 & 3.96 & 3.93\\\bottomrule
        \end{tabular}
      \end{scriptsize}
    \end{subtable}%
    \quad
    \begin{subtable}{0.5\linewidth}
      \centering
      \begin{scriptsize}
        \caption{$N=4$}
         \begin{tabular}{@{}lccc@{}}
         \toprule
         $K$ & $L^2(\phi)$ & $L^2(\phix)$ & $L^2(\phiy)$\\
         \midrule
          $4^2$ & 2.51E-04 & 8.81E-04 & 1.63E-03\\
          $8^2$ & 8.52E-06 & 2.88E-05 & 5.45E-05\\
         $16^2$ & 2.77E-07 & 9.12E-07 & 1.76E-06\\
         $32^2$ & 8.85E-09 & 2.85E-08 & 5.60E-08\\\midrule
         avg. EOC & 4.93 & 4.97 & 4.94\\\bottomrule
        \end{tabular}
      \end{scriptsize}
    \end{subtable}
\end{table}

\linelabel{lne:rk3star_speed_start}
\label{rev1:r3q3}\addrev{any}{To close the discussion regarding the high-order DG approximation of the hyperbolic diffusion equations, we present results to demonstrate how the RK3S* time integrator discussed in \secref{3STime} can accelerate convergence to steady state and approximate the solution of a Poisson problem. We compare the number of pseudotime steps for the RK3S* integrator against that required by the standard low-storage CK45 time integrator in \tabref{HypDiffCost}. A similar study of convergence acceleration via RK3S* is carried out in \secref{jeans} for a self-gravitating gas configuration. For this comparison we take $\CFLHG=1.0$ for the RK3S* integrator and $\CFLHG=0.5$ for the CK45 time integrator, which are the largest stable CFL numbers for both polynomial orders and each scheme, respectively. These CFL values produce identical convergence results (up to machine precision) but serve to illustrate the acceleration afforded by a RK integration technique specifically designed for time integration of the hyperbolic diffusion system.} 

\begin{table}[!htb]
   \centering
   \caption{\addrev{any}{Comparison of number of pseudotime steps needed for steady-state convergence
     of the hyperbolic diffusion equations to $\mathtt{tol}=10^{-10}$ for the manufactured solution
     test case \eqnref{PoissonManufac} for an increasing resolution of elements and fixed polynomial
     degree.}}
    \label{tab:HypDiffCost}
    \begin{subtable}{0.5\linewidth}
      \centering
        \begin{scriptsize}
        \caption{$N=3$}
         \begin{tabular}{@{}lccc@{}}
         \toprule
         $K$ & CK45 & RK3S* & Reduction \\
         \midrule
          $4^2$ & 793 & 397 & 49.9\%\\
          $8^2$ & 1587 & 794 & 49.9\%\\
         $16^2$ & 3180 & 1588 & 50.0\%\\
         $32^2$ & 6388 & 3185 & 50.1\%\\\bottomrule
        \end{tabular}
      \end{scriptsize}
    \end{subtable}%
    \begin{subtable}{0.5\linewidth}
      \centering
      \begin{scriptsize}
        \caption{$N=4$}
         \begin{tabular}{@{}lccc@{}}
         \toprule
         $K$ & CK45 & RK3S* & Reduction\\
         \midrule
          $4^2$ & 993 & 629 & 36.7\%\\
          $8^2$ & 1989 & 995 & 49.9\%\\
         $16^2$ & 4012 & 1999 & 50.2\%\\
         $32^2$ & 8338 & 4137 & 50.4\%\\\bottomrule
        \end{tabular}
      \end{scriptsize}
    \end{subtable}
\end{table}
\addrev{any}{We find that the optimized RK3S* time integration technique reduces the computational effort, measured with the number of pseudotime steps in the explicit Runge-Kutta solver, by approximately a factor of two for this problem setup. The wall clock time scales in tandem with the number of pseudotime steps required for convergence because the five-stage CK45 and RK3S* schemes have the same computational cost per step. This overall reduction of computational effort is due, in principal, to the increased CFL number (and subsequent explicit time step size) of the optimized RK3S* scheme as discussed in \secref{3STime}.}
\linelabel{lne:rk3star_speed_end}

\subsubsection{Coupled compressible Euler and gravity solver}\label{sec:coupleConv}

As a final verification test we apply the method of manufactured solutions to demonstrate the accuracy of the coupled DG simulation for compressible Euler with gravity. For completeness, we examine two coupling strategies:
\begin{itemize}
\setlength\itemsep{0em}
\item Updating the gravity system once in each RK \textit{stage} of the compressible Euler solver \addrev{2}{as described in \alref{gravEveryStage}},\linelabel{lne:new_alg_ref1}
\item Updating the gravity system once in each RK \textit{time step} of the compressible Euler solver \addrev{2}{as described in \alref{gravEveryStep}},\linelabel{lne:new_alg_ref2}
\end{itemize}
and show that their respective accuracies differ greatly.

Just as in \secref{euler_eoc}, we take the domain to be $\Omega=[0,2]^2$ with periodic boundary conditions, set $\gamma=2$ and take the manufactured solution for the compressible Euler variables to be \eqnref{EulerManufac}.
From the density solution ansatz in \eqnref{EulerManufac}, we take the manufactured solution of the gravitational potential to be
\begin{equation}\label{eqn:GravManufac}
\phi(x,y) = -\frac{2}{\pi}\,\frac{1}{10}\sin(\pi(x+y-t)) = -\frac{2}{\pi}(\rho - 2).
\end{equation}
This solution for $\phi(x,y)$ and its gradient 
\begin{equation}
\phix = \phiy = -\frac{1}{5}\cos(\pi(x+y-t)),
\end{equation}
are also periodic in the considered domain. Further, we compute
\begin{equation}\label{eq:gravPoisson}
-\vec{\nabla}^2\phi = -(\phi_{xx} + \phi_{yy}) = -4\pi(\rho - 2) = -4\pi\rho + 8\pi,
\end{equation}
which solves the gravitational Poisson problem \eqnref{Poisson} with the gravitational constant $G=1$ and a constant residual term of $8\pi$. It is straightforward to compute the remaining residual terms for the compressible Euler equations with gravity \eqnref{eulerPlusGravity} to be
\begin{equation}
 \vect{u}_t + \vec{\nabla} \cdot \bvect{f}(\vect{u}) =
\vect{s}(\vect{u})
+
\begin{bmatrix}
 \frac{1}{10} \pi \cos(\pi(x+y-t))\\[0.1cm]
 \frac{1}{10} \pi \cos(\pi(x+y-t))\\[0.1cm]
 \frac{1}{10} \pi \cos(\pi(x+y-t))\\[0.1cm]
\frac{1}{10} \pi \cos(\pi(x+y-t))\left[1 + \frac{2}{\pi}(2 + \frac{1}{10}\sin(\pi(x+y-t)))\right]\\[0.1cm]
\end{bmatrix},
\end{equation}
where $\vect{s}(\vect{u})$ are the source terms proportional to the gravity potential as in \eqnref{eulerPlusGravity}.

As in the two previous subsections we choose two polynomial orders $N=3$ and $N=4$ to demonstrate the accuracy of the DG solver for the couple simulations.
The  manufactured solution for the compressible Euler equations with gravity is run to $T=0.5$ as a final time. Note, the update of the gravitational potential and gradient variables must reach the prescribed tolerance in either every RK \textit{stage} or in every \textit{time step} of the compressible Euler solver, depending on the coupling strategy we select. 

In \tabref{couple_eoc_every_step}, we present the EOCs for the coupled manufactured solution test case where the gravity variables are updated in every RK \textit{stage} using the RK3S* scheme. The discrete $L^2$ errors for the conservative Euler variables as well as the hyperbolic diffusion variables $\phi$, $\phix$, and $\phiy$ are computed on uniform Cartesian meshes of increasing resolution. We see that this coupling strategy preserves the high-order accuracy of both DG solvers because the EOC of all solution variables is the optimal convergence order. We also ran the same convergence test configurations using CK45 to integrate both the hydrodynamic variables and hyperbolic gravity variables. The computed $L^2$ errors in all seven variables as well as their respective EOC were nearly identical to those obtained using the mixed CK45 for Euler and RK3S* for hyperbolic gravity. This confirms the validity of the RK3S* method derived in \secref{3STime}. Furthermore, this convergence test demonstrates that the choice of the explicit RK scheme used to drive the hyperbolic gravity system to steady state has no significant influence on the overall solution accuracy.
\begin{table}[!htb]
    \caption{Convergence test for manufactured compressible Euler with gravity solution \eqnref{EulerManufac} and \eqnref{GravManufac} coupled within every RK \textit{stage}. This coupling strategy retains the optimal convergence order for all variables.}
    \label{tab:couple_eoc_every_step}
    \begin{subtable}{\linewidth}
      \centering
      \begin{scriptsize}
        \caption{$N=3$}
         \begin{tabular}{@{}lccccccc@{}}
         \toprule
         $K$ & $L^2(\rho)$ & $L^2(\rho v_1)$ & $L^2(\rho v_2)$ & $L^2(E)$ & $L^2(\phi)$ & $L^2(\phix)$ & $L^2(\phiy)$\\
         \midrule
          $4^2$ & 4.37E-04 & 4.69E-04 & 4.69E-04 & 9.72E-04 & 1.64E-04 & 8.33E-04 & 8.33E-04\\
          $8^2$ & 2.43E-05 & 2.60E-05 & 2.60E-05 & 5.09E-05 & 9.90E-06 & 5.65E-05 & 5.65E-05\\
         $16^2$ & 1.06E-06 & 1.37E-06 & 1.37E-06 & 2.65E-06 & 6.63E-07 & 3.77E-06 & 3.77E-06\\
         $32^2$ & 4.73E-08 & 8.03E-08 & 8.03E-08 & 1.56E-07 & 4.33E-08 & 2.44E-07 & 2.44E-07\\\midrule
         avg. EOC & 4.39 & 4.17 & 4.17 & 4.20 & 3.96 & 3.91 & 3.91\\\bottomrule
         \end{tabular}
        \end{scriptsize}
    \end{subtable}%
    \vspace{0.25cm}
    \\
    \begin{subtable}{\linewidth}
      \centering
      \begin{scriptsize}
        \caption{$N=4$}
         \begin{tabular}{@{}lccccccc@{}}
         \toprule
         $K$ & $L^2(\rho)$ & $L^2(\rho v_1)$ & $L^2(\rho v_2)$ & $L^2(E)$ & $L^2(\phi)$ & $L^2(\phix)$ & $L^2(\phiy)$\\
         \midrule
          $4^2$ & 3.50E-05 & 3.38E-05 & 3.38E-05 & 6.59E-05 & 1.15E-05 & 6.31E-05 & 6.31E-05\\
          $8^2$ & 7.99E-07 & 9.00E-07 & 9.00E-07 & 1.71E-06 & 3.74E-07 & 2.11E-06 & 2.11E-06\\
         $16^2$ & 1.95E-08 & 2.49E-08 & 2.49E-08 & 4.78E-08 & 1.23E-08 & 6.95E-08 & 6.95E-08\\
         $32^2$ & 5.31E-10 & 7.73E-10 & 7.73E-10 & 1.44E-09 & 4.03E-10 & 2.25E-09 & 2.25E-09\\\midrule
         avg. EOC & 5.34 & 5.14 & 5.14 & 5.16 & 4.93 & 4.93 & 4.93\\\bottomrule
         \end{tabular}
        \end{scriptsize}
    \end{subtable}
\end{table}

Next, we provide the EOCs in \tabref{N3eoc_every_stage} for polynomial order $N=3$ of the coupled
test case where the gravity variables are updated in every RK \textit{time step} \addrev{2}{of the compressible Euler simulation}. The results for
$N=4$ are omitted for brevity, but the results are similar. We see that the convergence order in
\tabref{N3eoc_every_stage} has dropped to first order for this coupling strategy. The error is no
longer dominated by spatial errors (as all previous results) but instead by the temporal discretization. This is because the gravitational potential is treated as ``fixed'' for the given RK stages before it is updated again. Thus, to improve the approximation accuracy requires one to shrink ${\CFLE}$ to mitigate the error introduced when treating the gravitational potential in this frozen way.
\begin{table}[!htp]
\caption{Convergence test for compressible Euler and gravity manufactured solution test case coupled in every RK \textit{time step} with polynomial order $N=3$. This demonstrates such a coupling technique introduces a first-order error into the approximation.}
\label{tab:N3eoc_every_stage}
\centering
\begin{scriptsize}
\begin{tabular}{@{}lccccccc@{}}
\toprule
$K$ & $L^2(\rho)$ & $L^2(\rho v_1)$ & $L^2(\rho v_2)$ & $L^2(E)$ & $L^2(\phi)$ & $L^2(\phix)$ & $L^2(\phiy)$\\
\midrule
  $4^2$ & 7.52E-03 & 7.55E-03 & 7.55E-03 & 1.68E-02 & 4.79E-03 & 1.51E-02 & 1.51E-02\\
  $8^2$ & 3.85E-03 & 3.90E-03 & 3.90E-03 & 8.67E-03 & 2.45E-03 & 7.69E-03 & 7.69E-03\\
$16^2$ & 1.95E-03 & 1.99E-03 & 1.99E-03 & 4.43E-03 & 1.24E-03 & 3.90E-03 & 3.90E-03\\
$32^2$ & 9.82E-04 & 1.01E-03 & 1.01E-03 & 2.24E-03 & 6.25E-04 & 1.96E-03 & 1.96E-03\\\midrule
avg. EOC & 0.98 & 0.97 & 0.97 & 0.97 & 0.98 & 0.98 & 0.98\\\bottomrule
\end{tabular}
\end{scriptsize}
\end{table}

To conclude, coupling of compressible Euler and hyperbolic gravity solvers within every RK \textit{stage} of the Euler solver is preferred, because it preserves the high-order accuracy of the DG spatial approximation. We will investigate further the influence of the coupling strategy on solution quality for a more practical example in self-gravitating flows in \secref{jeans}.

\subsection{Applications for self-gravitating gas dynamics}\label{sec:num_appl}

Beyond the verification test cases, we demonstrate the multi-physics capabilities of \trixi in simulating two self-gravitating flows. First, in \secref{jeans}, we consider the Jeans instability \cite{Jeans1902} that models perturbations and interactions between a gas cloud and gravity. In \secref{sedov}, we exercise the shock capturing and AMR capabilities of \trixi to simulate a self-gravitating variant of the Sedov blast wave.
\linelabel{lne:tolerance_start}
For these simulations we set the steady-state tolerance for the gravity solver to be $\mathtt{tol}=10^{-4}$ as is often done in astrophysical simulations \cite{FLASHug,ricker2008direct,huang1999fast}. 
\label{rev1:r2q4part1}\addrev{2}{Experience across self-gravitating applications and simulation codes reinforces that this choice of error tolerance offers a good balance between solution quality and the performance of the gravity solver \cite{Johannes}.}
\addrev{2}{We note, however, that determining a suitable tolerance value for arbitrary simulation setups is a non-trivial
task and presents a challenge that is similar to finding optimal parameters for classical
iterative solvers or implicit time integration schemes. Furthermore,}
\linelabel{lne:tolerance_end}
unless stated otherwise, the coupling of the two solvers is performed in every RK \textit{stage} to preserve the high-order spatial accuracy of the DG approximations.
All performance tests were conducted on a machine with an Intel Core i7-6850K CPU at 3.60~GHz and
32~GiB main memory. The presented numbers represent the minimum run times out of four separate
measurements and were obtained by executing Julia on one thread with bounds checking disabled.

\subsubsection{Jeans gravitational instability}\label{sec:jeans}

A simple example for an instability in a self-gravitating, thermally supported
interstellar cloud was first described by Jeans \cite{Jeans1902}. The linear instability mode is particularly useful to test the coupling of hydrodynamics to gravity, since it is one of the
few problems with periodic gravitational potential for which there exists an analytical solution for comparison, e.g., \cite{FLASHug,derigs2016novel}. Approximating the Jeans instability allows for the (numerical) study of pressure-dominated and gravity-dominated flows as well as the oscillation of the self-gravitating gas cloud between the two limits.

The domain is $\Omega=[0,1]^2$ with periodic boundary conditions for the hydrodynamics as well as gravity components and $\gamma=5/3$. We summarize the initial conditions in Centimeter-Gram-Seconds (CGS) units here but further details on their interpretation and derivation can be found in \cite{hubber2006resolution,FLASHug}. 

Consider a static medium, with uniform density $\rho_0$ and pressure $p_0$ at rest. Further, assume that any fluctuations between density and pressure occur adiabatically, such that $p_0 = \gamma \rho_0$. Now, we suppose that this uniform medium is initially perturbed so that
\begin{equation}\label{eqn:JeansPerturb}
\begin{aligned}
\rho &= \rho_0 + \rho_1=\rho_0\left[1 + \delta_0\cos(\vec{k}\cdot\vec{x})\right],\\[0.1cm]
p &= p_0 + p_1 = p_0\left[1 + \delta_0\gamma\cos(\vec{k}\cdot\vec{x})\right],\\[0.1cm]
\vec{v} &= \vec{0},
\end{aligned}
\end{equation}
where $\delta_0 = 10^{-3}$ is the amplitude of the perturbation and $\vec{k} = (4\pi,0)^\transp$ [cm$^{-1}$] is the wave vector that dictates the perturbation mode with the associated wave number $k^2 = \vec{k} \cdot\vec{k}$. The background medium values are taken to be $\rho_0 = 1.5 \cdot 10^7$ [g cm$^{-3}$] and $p_0 = 1.5 \cdot 10^7$ [dyn cm$^{-2}$]. The gravitational potential due to the perturbed density is given by \cite{binney2011galactic}
\begin{equation}
-\vec{\nabla}^2\phi = -4\pi G(\rho - \rho_0) = -4\pi G\rho_1,
\end{equation}
with $G = 6.674\cdot10^{-8}$ [cm$^3$g$^{-1}$s$^{-2}$] as the gravitational constant. For the initialization of the hyperbolic gravity solver, we assume constant state for the gravitational potential of $\phi = \delta_0\rho_0$ and a constant zero state for the auxiliary gradient variables in \eqnref{hypGrav}.

It is possible to obtain the dispersion relation of the self-gravitating fluid perturbation \eqnref{JeansPerturb} by examining a plane wave solution in Fourier space \cite{hubber2006resolution} to find
\begin{equation}\label{eqn:JeansDispersion}
\omega^2 = c_0^2k^2 - 4\pi G \rho_0,
\end{equation}
where $c_0 = \sqrt{\gamma p_0/\rho_0}$ [cm s$^{-1}$] is the ambient sound speed. From \eqnref{JeansDispersion} we define the Jeans wave number
\begin{equation}
k_J = \frac{\sqrt{4\pi G \rho_0}}{c_0} \approx 2.75,
\end{equation}
for the considered initial value configuration. The Jeans wave number is of critical importance because it separates between two physically relevant regimes. When $k>k_J$ the perturbation varies periodically in time and the equilibrium is \textit{stable}. That is, the perturbation amplitude simply oscillates transferring energy into gravitational potential (and vice versa). It does not increase with time. If, however, $k<k_J$ the perturbation is unstable and the amplitude grows exponentially in time, leading to a gas cloud that becomes denser and denser, eventually resulting in gravitational collapse \cite{bonnor1957jeans,Chandrasekhar1961}. For the perturbation parameters in \eqnref{JeansPerturb}, $k = 4\pi > k_J$ and the resulting perturbation is stable.

We simulate the Jeans gravitational perturbation \eqnref{JeansPerturb} with the \trixi multi-physics solver on a uniform $16\times 16$ Cartesian mesh with polynomial order $N=3$ in each spatial direction, resulting in $64^2$ degrees of freedom for each equation variable. We run the simulation up to a final time of $T=5.0$ [s], corresponding to approximately sixteen full oscillations of the perturbation. The standard DGSEM approximation is used for both the compressible Euler solver and the hyperbolic gravity solver.

To examine the behaviour of the Jeans instability, we investigate the bulk values of the
kinetic, internal, and potential energies defined by
\begin{equation}\label{eqn:threeEnergies}
  E_{\texttt{kin}} = \int\frac{\rho}{2}(v_1^2 + v_2^2)\diff \Omega, \qquad
  E_{\texttt{int}} = \int\frac{p}{\gamma -1}\diff \Omega, \qquad
  E_{\texttt{pot}} = \int\rho\phi\diff \Omega.
\end{equation}
We integrate the three energies \eqnref{threeEnergies} over the entire domain by applying LGL
quadrature. Then, we compare these approximate bulk energies against available analytical profiles
\cite{FLASHug,derigs2016novel}. 

\begin{figure}[H]
\centering
\includegraphics[width=0.79\textwidth]{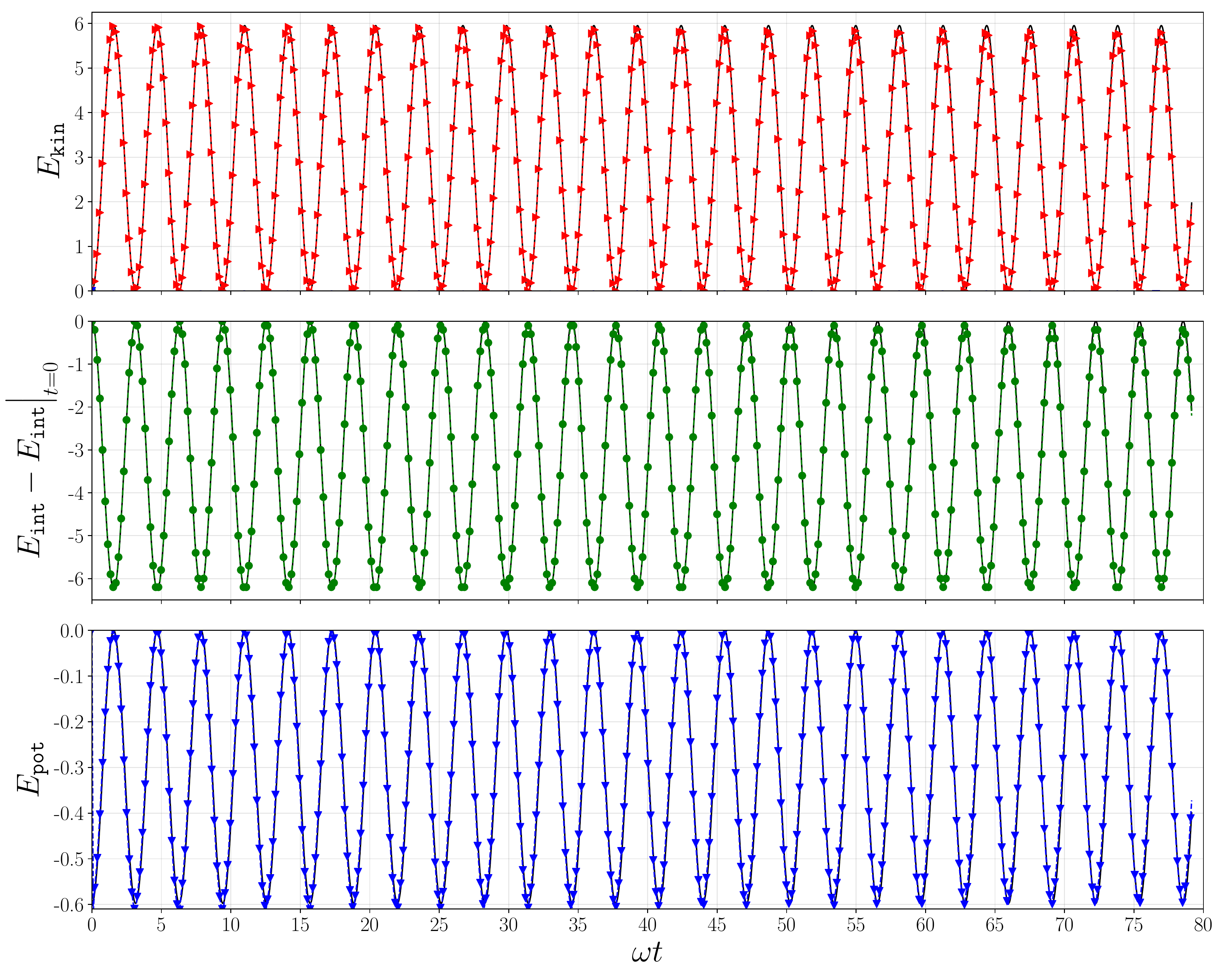}
\caption{Evolution of the computed kinetic (\protect\includegraphics[scale=0.2]{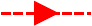}), internal (\protect\includegraphics[scale=0.2]{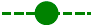}), and potential (\protect\includegraphics[scale=0.2]{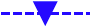}) energies for the Jeans instability using polynomial order $N=3$ on a uniform $16\times 16$ mesh. The analytical (\protect\includegraphics[scale=0.225]{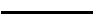}) energies are included for reference. Multi-physics coupling is done in every RK \textit{stage} of the compressible Euler solver with $\CFLE=0.5$ and $\CFLHG = 0.8$ using CK45.}
\label{fig:jeans_every_stage}
\end{figure}

In the first simulation of the Jeans gravitational perturbation we use the CK45 time integration scheme for the compressible Euler as well as hyperbolic gravity solvers. We select the explicit time step for the compressible Euler solver with $\CFLE=0.5$ and $\CFLHG = 0.8$ for the hyperbolic gravity solver. In \figref{jeans_every_stage} we show the resulting kinetic, internal and potential energy profiles as functions of $\omega t$. We see that the multi-physics approximation captures the amplitude and phase of the oscillatory solution very accurately over time, as one expects due to the excellent dissipation and dispersion properties of the DGSEM \cite{Airnsworth2004}.

The Jeans instability problem offers an interesting middle ground to investigate computational efficiency and solution accuracy. This is because the Jeans instability is a more physically relevant test setup than the manufactured solution from \secref{coupleConv}, but still possesses analytical energy profiles that we can compare against. As such, we run another simulation that employs the alternative gravity coupling procedure of ``freezing'' the gravitational potential within each RK \textit{time step} and evolving the hydrodynamic quantities. We demonstrated in \secref{coupleConv} that this introduces a first-order temporal error into the approximation; however, it is interesting to examine how such coupling influences the solution quality for the Jeans instability test case.

\begin{figure}[H]
\centering
\includegraphics[width=0.79\textwidth]{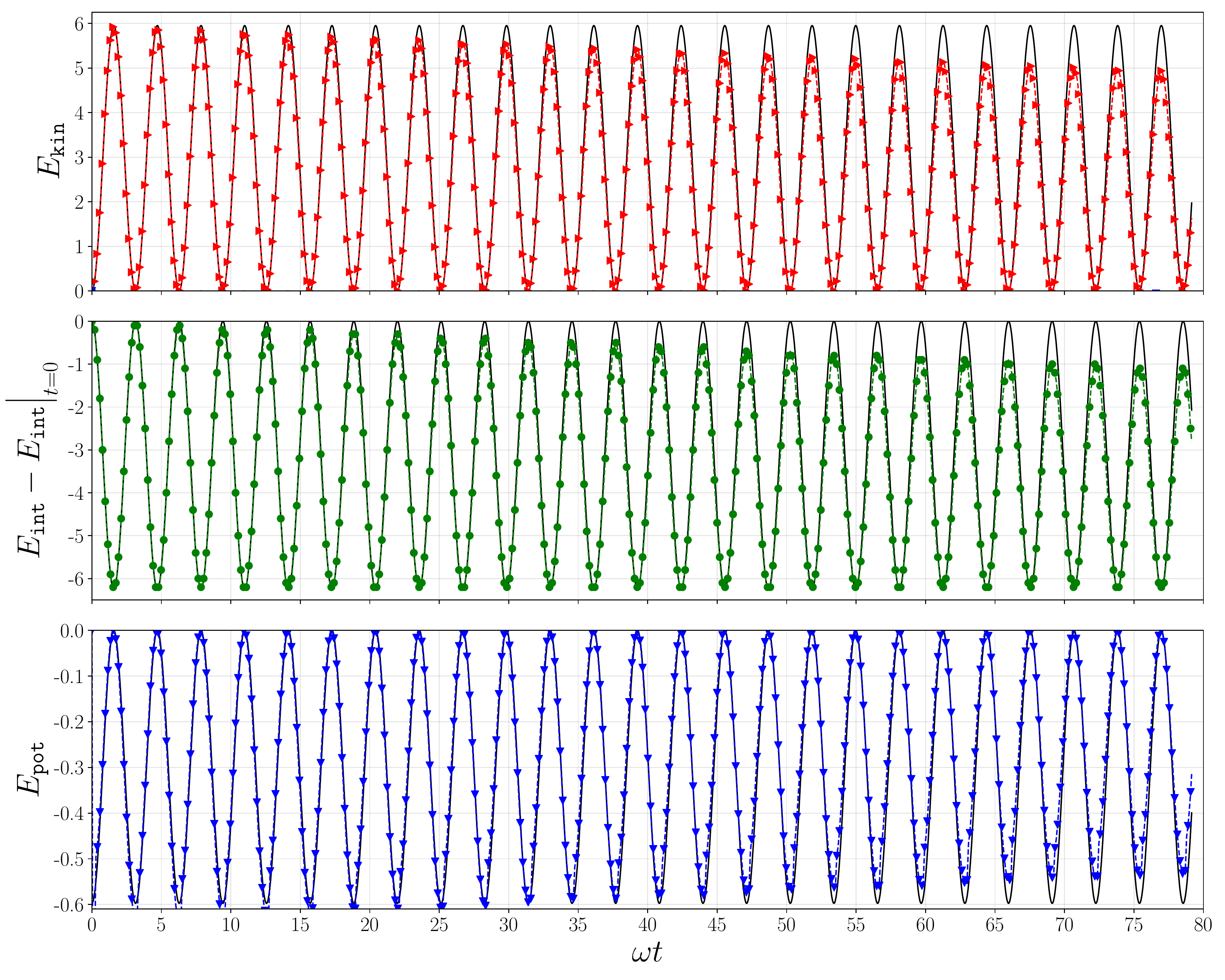}
\caption{Evolution of the kinetic (\protect\includegraphics[scale=0.2]{figures/jeans/label_ekin.png}), internal (\protect\includegraphics[scale=0.2]{figures/jeans/label_eint.png}), and potential (\protect\includegraphics[scale=0.2]{figures/jeans/label_epot.png}) energies for the Jeans instability using polynomial order $N=3$ on a uniform $16\times 16$ mesh. The analytical (\protect\includegraphics[scale=0.225]{figures/jeans/label_exact.png}) energies are included for reference. Multi-physics coupling is done in every RK \textit{time step} of the compressible Euler solver with $\CFLE=0.5$ and $\CFLHG = 0.8$ using CK45. There is a loss in energy amplitudes due to the first-order coupling.}
\label{fig:jeans_every_step}
\end{figure}

As such, we again run the multi-physics solver using CK45 time integration for the compressible Euler and hyperbolic gravity solvers with $\CFLE=0.5$ and $\CFLHG = 0.8$ for time step selection. This time, however, we update the gravitational potential and its gradients after every RK \textit{time step} of the Euler solver. We present the evolution of the kinetic, internal and potential energy profiles in \figref{jeans_every_step} as functions of $\omega t$. We, again, see that the dispersion errors are very small for the DG approximation of the Jeans instability with this alternative coupling. But, as time progresses, there is a noticeable loss in amplitude of the different energies due to the first-order errors introduced into the approximation. Such errors can be removed by taking a very small value of $\CFLE$ \cite{derigs2016novel} (for example selecting $\CFLE=0.01$ produces results nearly identical to \figref{jeans_every_stage}), effectively reducing the temporal errors to be of lower magnitude to reveal dominate spatial behavior. However, such a small value of $\CFLE$ is computationally prohibitive and further reinforces our finding from \secref{coupleConv} that a multi-physics coupling procedure that preserves the underlying order of the spatial approximation is also desirable for practical simulations.

We also investigate the effect that the optimized method RK3S* from \secref{3STime}
has in reducing the computational effort for the hyperbolic gravity
solver. 
\linelabel{lne:grav_cycle_start}
\label{rev1:r2q2}\addrev{2}{To do so, we introduce a shorthand cost measurement deemed a gravity \textit{sub-cycle}. This cost measure corresponds to the assembly and evolution of the hyperbolic gravity system by one complete time step in pseudotime within the gravity update loop (see \figref{trixi_flowchart}). We note that each gravity sub-cycle consists of five RK stages for either the CK45 or RK3S*  time integration schemes.}
Recall, in each RK \textit{stage} \addrev{2}{of the compressible Euler simulation} the hyperbolic gravity solver must
\addrev{2}{evolve in pseudotime} until the steady-state tolerance of $\mathtt{tol}=10^{-4}$ is reached. 
\linelabel{lne:pseudotime_4}
Thus, an obvious first attempt to improve the computational efficiency of the multi-physics implementation is to reduce the number of \addrev{2}{gravity} 
\linelabel{lne:grav_cycle_end}
sub-cycles.
We again simulate the Jeans instability with $\CFLE=0.5$ but compare the
sub-cycle counts for the hyperbolic gravity solver with the CK45 scheme with
$\CFLHG=0.8$ against RK3S* with $\CFLHG=1.2$. These CFL values led to the largest stable explicit time steps that still produced meaningful simulation results. 

\begin{figure}[H]
     \centering
     \begin{subfigure}[b]{0.49\textwidth}
         \centering
         {\includegraphics[width=\textwidth]{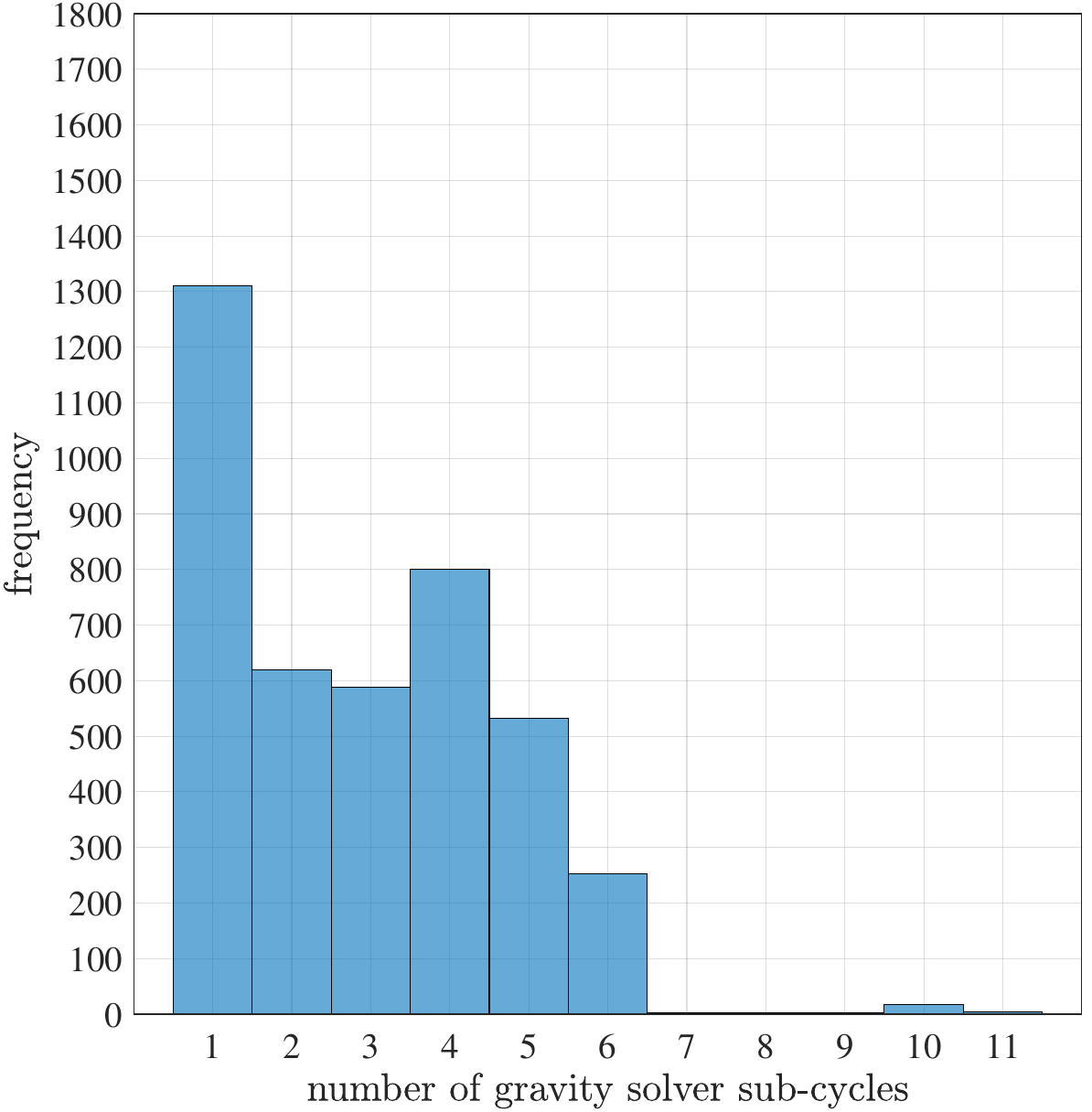}}
         \caption{CK45, $\CFLE=0.5$, $\CFLHG=0.8$}
         \label{fig:jeans_hist_rk45}
     \end{subfigure}
     \hfill
     \begin{subfigure}[b]{0.49\textwidth}
         \centering
         {\includegraphics[width=\textwidth]{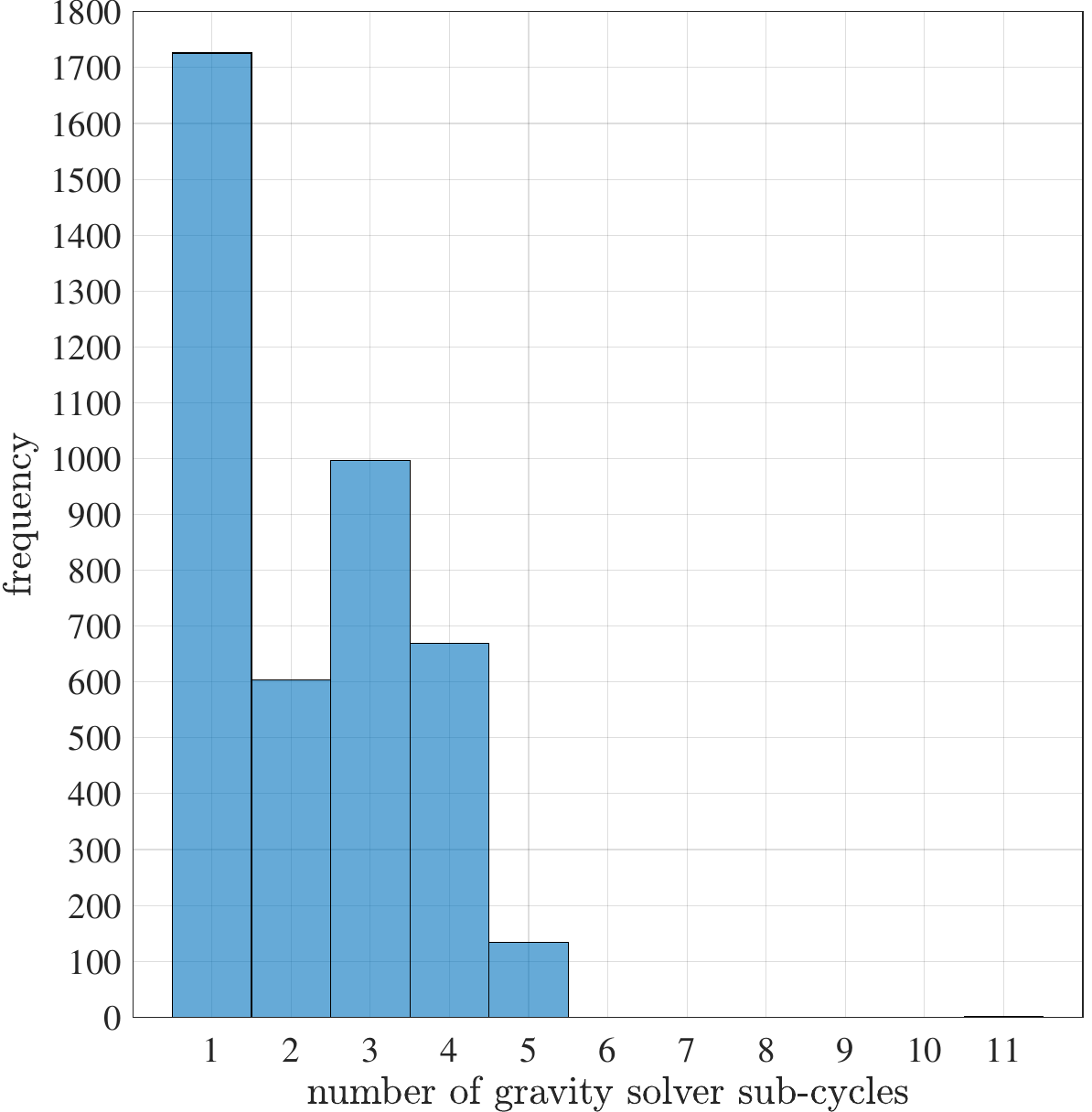}}
         \caption{RK3S*, $\CFLE=0.5$, $\CFLHG=1.2$}
         \label{fig:jeans_hist_3Sstar}
     \end{subfigure}
     \caption{Two runs of the Jeans instability test case with polynomial order $N=3$ on a uniform $16\times 16$ mesh with different explicit time integration methods. Histograms present the frequency of different \addrev{2}{gravity} sub-cycles counts needed for the hyperbolic gravity solver to reach the steady-state tolerance $\mathtt{tol}=10^{-4}$ in each RK \textit{stage}. Comparison can be made between the number and distribution of gravity solver sub-cycles necessary for CK45 (left) or RK3S* (right) optimized for the hyperbolic diffusion system \eqnref{hypGrav}.}
     \label{fig:jeansHist}
\end{figure}

Histograms in
\figref{jeansHist} visualize the sub-cycle frequency of the hyperbolic gravity
solver for these two time integration techniques. Note we have removed
outlier \addrev{2}{gravity} sub-cycle values that occur only {once} (\eg, for the initial solve in
the first time step) to better illustrate the trend and overall effort of the
hyperbolic gravity solver runs.

The results shown in \figref{jeansHist} demonstrate that the number of gravity solver sub-cycles is concentrated near one or two iterations for either time integration scheme. So, using the hyperbolic gravity variables from the previous compressible Euler RK stage as the initial guess for the next hyperbolic gravity solve works well to keep the number sub-cycles to update the gravitational potential small. However, the spread of the sub-cycle iteration number is wider for the standard CK45 time integrator compared to RK3S* that was optimized to take larger explicit time steps for the hyperbolic gravity problem. Apart from the concentration of the sub-cycle distribution it is also noteworthy that the raw number of hyperbolic gravity sub-cycles for CK45 was 12,055 whereas RK3S* required only 9,344. Therefore, not only is the distribution of \addrev{2}{gravity} sub-cycles more skewed toward lower values for RK3S*, but the computational effort for the gravity solver is also decreased by approximately \SI{22}{\percent}.
We would like to emphasize that this performance improvement does not negatively impact solution
accuracy, as discussed in \secref{coupleConv}, and that it can be applied immediately if 3S*
RK schemes are already supported by the implementation.

\subsubsection{Sedov explosion with self-gravity}\label{sec:sedov}

As a final numerical example we consider a modification of the Sedov blast wave problem that incorporates the effects of gravitational acceleration \cite{FLASHug}. The hydrodynamic setup of the Sedov explosion \cite{sedov1993similarity} is a difficult one, as it involves strong shocks and complex fluid interactions. We include it to demonstrate the shock capturing and AMR capabilities of \trixi to resolve the cylindrical Sedov blast wave. Additionally, it highlights that the treatment of the gravitational potential as a hyperbolic system \eqnref{hypGrav} is immediately amenable to AMR through a standard mortar method. No further considerations are necessary to approximate the gravitational potential on non-conforming meshes.

The initial configuration of the Sedov problem deposits the explosion energy $E$ into a single point in a medium of uniform ambient density $\rho_{\mathrm{am}}$ and pressure $p_{\mathrm{am}}$. In practice, the initialization of the Sedov problem is delicate because this energetic area is typically smaller than the grid resolution. Therefore, we follow an approach similar to Fryxell et al.~\cite{fryxell2000flash} to convert the explosion energy into a pressure contained within a resolvable area center of radius $r_{\mathrm{ini}}$ by
\begin{equation}\label{eqn:Sedov_p_ini}
p_{\mathrm{ini}} = \frac{(\gamma - 1)E}{\pi r_{\mathrm{ini}}}.
\end{equation}
This pressure is then used for the discretization points where $r < r_{\mathrm{ini}}$. For the
simulation we choose $r_{\mathrm{ini}}$ to be four times as large as the initial grid spacing, which helps to minimize effects due to the Cartesian geometry of the computational grid.

We consider the Sedov explosion parameters in CGS units to be $p_{\mathrm{am}} = 10^{-5}$ [dyn
cm$^{-2}$], $E=1$ [erg], and $v_1 = v_2 = 0$ [cm s$^{-1}$], and we set $\gamma = 1.4$. The
computational domain $\Omega = [-4,4]^2$ [cm$^2$] is discretized by an adaptive mesh with a minimum
element length $h = 0.03125$ [cm]. Thus, the value
$r_{\mathrm{ini}} = 0.125$ [cm] is used in \eqnref{Sedov_p_ini}. The gravitational constant is $G
= 6.674\cdot10^{-8}$ [cm$^3$g$^{-1}$s$^{-2}$]. For the gravitational potential it is customary to
assume that it vanishes at large distances away from a localized region of non-zero density, e.g.,
\cite{hubber2018gandalf,FLASHug,katz2016white}. Therefore, we localize the ambient density,
$\rho_{\mathrm{am}}$, to be contained in a disc of radius $r_{\rho} = 1$ [cm] such that
\begin{equation}\label{eqn:local_density}
\rho_{\mathrm{am}}\;[\text{g cm}^{-3}]= \begin{cases}
1, & r \leq r_{\rho}\\
10^{-5} , & r > r_{\rho}. \\
\end{cases}
\end{equation}
The transition between inner and ambient state for both pressure and density is smoothened by a
logistic function with steepness $k = 150$ [cm$^{-1}$].
We set Dirichlet boundary conditions at the four edges of the domain to be the
ambient flow states for the hydrodynamic variables and zero for the hyperbolic gravity system. 

\begin{figure}[!htp]
  \centering
  \includegraphics[width=\textwidth]{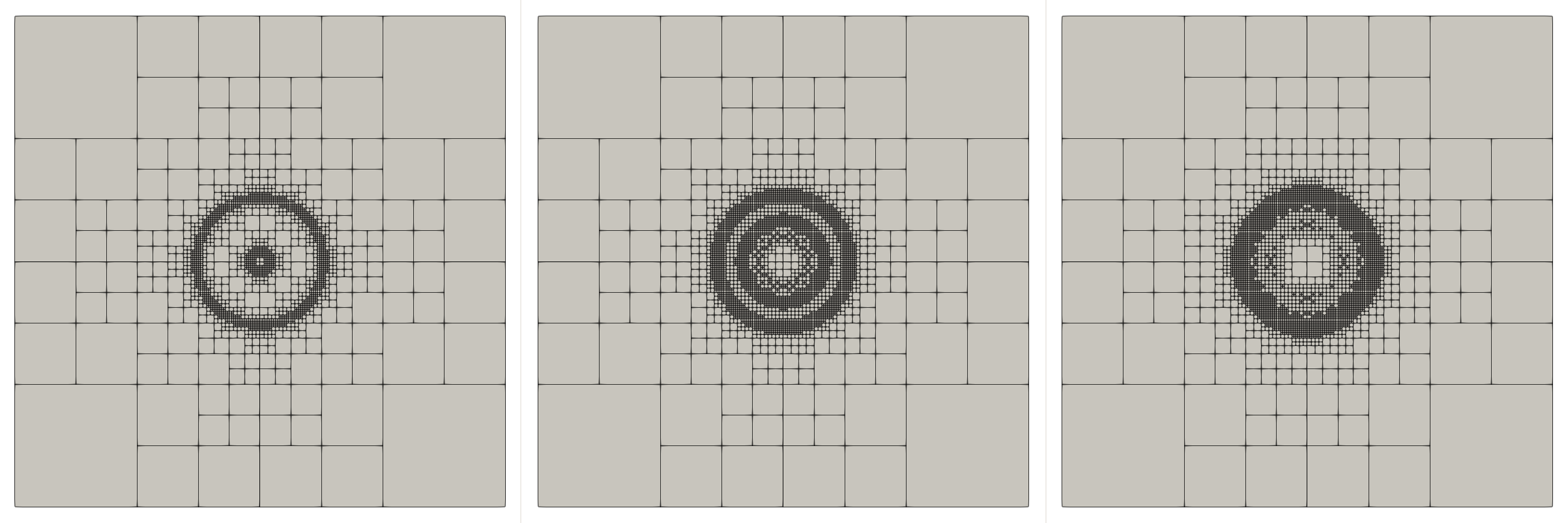}
  \caption{Adaptive mesh for self-gravitating Sedov explosion with seven different refinement levels
  at $T=0.0$ (left), $T=0.5$ (center), and $T=1.0$ (right).}
  \label{fig:sedov_amr_meshes}
\end{figure}

We approximate the hydrodynamic and hyperbolic gravity solutions with polynomial order $N=3$ and run
the simulation to a final time of $T=1.0$. We select a time step with $\CFLE=0.5$ for CK45 and
$\CFLHG=1.2$ using RK3S*. The compressible Euler solver uses the split-form DGSEM feature of \trixi
with the shock capturing scheme by Hennemann et al.\
\cite{HENNEMANN2021109935} as described in \secref{numerical_methods}, while the hyperbolic
gravity solver uses the standard DGSEM formulation.
During the simulation, the mesh is adaptively coarsened and refined after every time
step of the compressible Euler solver. The mesh resolution spans seven refinement levels, where $h =
2.0$~[cm] at the base level ($l = 2$) and $h = 0.03125$~[cm] at the highest refinement level ($l =
8$), see \figref{sedov_amr_meshes}. We evaluate the same indicator function for AMR as for shock
capturing.  Elements with $\alpha_\text{AMR} > 0.0003$ are assigned a target refinement level of $l
= 8$, all other elements are assigned a target level of $l = 2$. The value of the AMR indicator
$\lambda$ is based on whether an element matches its target level or needs to be adapted, and the
2:1 balancing algorithm automatically ensures a smoothly varying mesh resolution.  For details, see
the description in \appref{algorithms_implementation}.

\begin{figure}[!htp]
     \centering
     \begin{subfigure}[!htp]{0.49\textwidth}
         \centering
         \includegraphics[width=\textwidth]{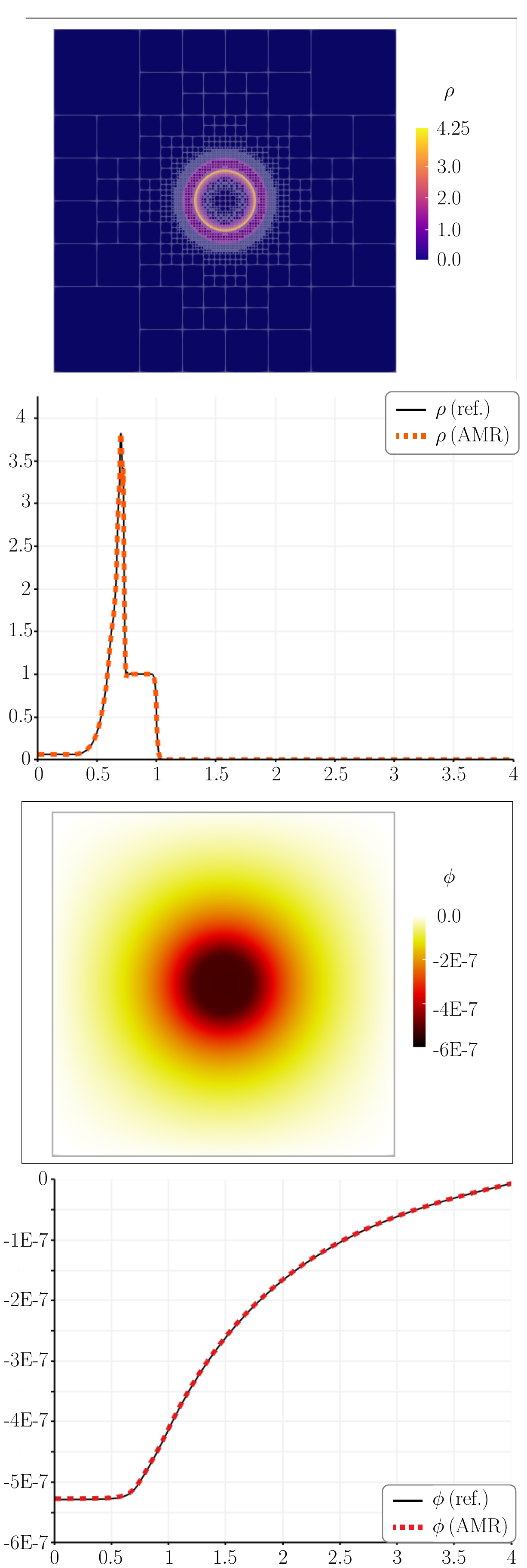}
         \caption{$T = 0.5$}
         \label{fig:sedov_t_half}
     \end{subfigure}
     \hfill
     \begin{subfigure}[!htp]{0.49\textwidth}
         \centering
         \includegraphics[width=\textwidth]{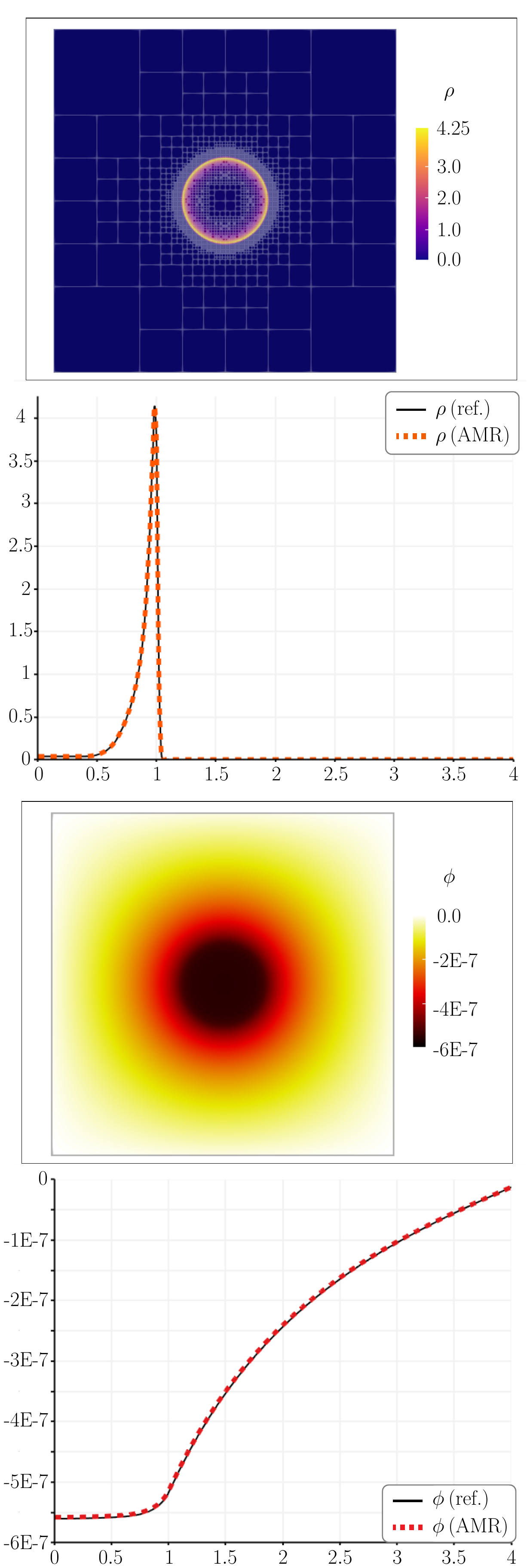}
         \caption{$T = 1.0$}
         \label{fig:sedov_t_one}
     \end{subfigure}
     \caption{Self-gravitating Sedov explosion approximated with polynomial order $N=3$ on an AMR
     grid with four levels of possible refinement. (first and third row) Two dimensional plots of
   the density and gravitational potential at two times. The white overlay of squares in the density
 plots shows the AMR grids used by both solvers. (second and last row) One dimensional slices of the solutions along a line from the origin to the edge of the domain in the positive $x$-direction.}
     \label{fig:sedovs}
\end{figure}

We present the approximate density and gravitational potential at the intermediate time $T=0.5$ and final time $T=1.0$ in the first and third rows of \figref{sedovs}. The density pseudocolor plots, in the top row, also give an overlay of the joint AMR grid that is shared by both solvers. The second and fourth rows of \figref{sedovs} extract a slice of the density and gravitational potential solutions along the horizontal line from the origin to the edge of the domain in the positive $x$-direction. The shapes of these one-dimensional profiles match well with the results for a similar self-gravitating Sedov explosion test available in \cite{FLASHug}. Furthermore, we create a reference solution for this test problem using a high-resolution run on a uniform mesh at the finest AMR level with polynomial order $N=3$ in each spatial direction. The one-dimensional profile slices in \figref{sedovs} reveal that the solutions on the adaptive mesh are virtually indistinguishable from the reference result.

As a final result, the computational effort of \trixi to solve the self-gravitating Sedov explosion
with a uniform and an adaptive mesh is presented in \tabref{sedov_timings}.
It is worth noting that the percentage figures of the gravity solver are comparable to other
multi-physics solvers. \label{rev2:r3q1}\linelabel{lne:grav_percent_again_start}\addrevtwo{3}{For
instance, the percentage of the total runtime spent in the gravity solver was reported for
other astrophysical simulation frameworks to be approximately
88--92\% in \cite[Table~1]{SPRINGEL200179}, 
43--78\% in \cite[Fig.~2]{Almgren_2010}, 
25--36\% in \cite[Fig.~C1]{VogelsbergerSijackiEtAl12}, 
85\% in \cite[Fig.~6]{HiraiNagakuraEtAl16}, 
34--72\% in \cite[Table~5]{WuenschWalchEtAl18}, or 
14--50\% in \cite[Fig.~6]{Moon_2019}. 
This underpins that the presented hyperbolic formulation of the gravity equation offers an
alternative solution strategy to existing gravity solvers.}
\linelabel{lne:grav_percent_again_end}
We further see that activating
AMR for this test case decreases the overall run time by approximately a factor of 19,
while there is minimal difference between the reference and adaptive solutions for the
self-gravitating Sedov blast wave. This demonstrates a novel advantage of non-conforming DG with AMR
capabilities and the hyperbolic gravity formulation. It is possible to accurately and rapidly
approximate the solution of an elliptic problem in a hyperbolic fashion.
\begin{table}[!htp]
\caption{Run time measurements for simulating the self-gravitating Sedov blast with a uniform
and an adaptive mesh. The uniform mesh is at refinement level $l = 8$, while the adaptive mesh is at
multiple levels $l \in [2, 8]$ with adaptation performed after every time step of the compressible
Euler solver.}
\label{tab:sedov_timings}
\centering
\begin{scriptsize}
\begin{tabular}{@{}lrrrr@{}}
\toprule
~ & \multicolumn{2}{c}{Uniform mesh} & \multicolumn{2}{c}{Adaptive mesh} \\
\cmidrule{1-1}\cmidrule(lr){2-3}\cmidrule(lr){4-5}
Compressible Euler solver   & $652.8 \,\textrm{s}$ & $46.3\%$ & $ 39.2 \,\textrm{s}$ & $52.4\%$ \\
Gravity solver              & $758.0 \,\textrm{s}$ & $53.7\%$ & $ 26.5 \,\textrm{s}$ & $35.4\%$ \\
AMR                         &                   -- &       -- & $  9.1 \,\textrm{s}$ & $12.1\%$ \\
\cmidrule{1-1}\cmidrule(lr){2-3}\cmidrule(lr){4-5}
Total                       &$1410.8 \,\textrm{s}$ & $100.0\%$& $ 74.8 \,\textrm{s}$ & $100.0\%$\\
\bottomrule
\end{tabular}
\end{scriptsize}
\end{table}

\section{Conclusions}
\label{sec:conclusions}

In this paper, we adapted Nishikawa's idea for solving the Poisson equation as a hyperbolic
diffusion system in the context of self-gravitating flows. 
\linelabel{lne:concl_motive_1_start}
\addrev{any}{In doing so, we converted 
the equation for a Newtonian gravitational potential into a system of hyperbolic gravity equations.}
\linelabel{lne:concl_motive_1_end}

Then, the hydrodynamics and hyperbolic gravity approximations  
were both built from a nodal discontinuous Galerkin method. This allowed us to treat the elliptic
problem for Newtonian gravity in a fully explicit hyperbolic fashion, and we verified that we are
able to obtain high-order accurate solutions for both single- and multi-physics simulations. By treating the elliptic
gravity problem in a purely hyperbolic way, its numerical approximation
inherits the full functionality and accuracy of the nodal discontinuous Galerkin
method. Most notably, the ability to do non-conforming approximations via mortars and to solve the
elliptic gravity problem on an adaptive grid is retained without any modifications. In addition,
coupling the flow and the gravity solver in \emph{each} Runge-Kutta stage preserves and carries over
the high-order benefits from the single-physics into the multi-physics context. Borrowing ideas from
the time integration community, we optimized a Runge-Kutta scheme to reduce computational
effort and to accelerate the hyperbolic gravity solution to steady state.

\linelabel{lne:concl_motive_2_start}
\addrev{any}{From the numerical results, we found that the fractions of computational effort for the compressible 
Euler and hyperbolic gravity solvers were comparable to other multi-physics solvers. 
Thus, for problems in self-gravitating gas dynamics, the use of hyperbolic diffusion methodology provides an interesting 
alternative to a classical iterative elliptic solver. The further development and optimization of the hyperbolic gravity multi-physics 
approach is the subject of ongoing research.}
\linelabel{lne:concl_motive_2_end}
All results in this paper
were obtained with the open-source simulation framework \trixi, which is available as a registered
Julia package. To allow reproducing our findings and to invite further collaboration, we also made
the corresponding setup files publicly available online \cite{SchlottkeLakemperWintersEtAl20}.

\section*{Acknowledgments}

Michael Schlottke-Lakemper thanks Stefanie Walch for valuable insights on the application and
performance of astrophysical simulation methods.
Andrew Winters thanks Johannes Markert for helpful discussions on gravitational solvers and astrophysical applications.
Gregor Gassner thanks the European Research Council for funding through the ERC Starting Grant ``An Exascale aware and Un-crashable Space-Time-Adaptive Discontinuous Spectral Element Solver for Non-Linear Conservation Laws'' (\texttt{Extreme}).
Andrew Winters was supported by Vetenskapsr{\aa}det, Sweden award number 2020-03642 VR.
Research reported in this publication was supported by the
King Abdullah University of Science and Technology (KAUST).
Funded by the Deutsche Forschungsgemeinschaft (DFG, German Research Foundation)
under Germany's Excellence Strategy EXC 2044-390685587, Mathematics Münster:
Dynamics-Geometry-Structure.

\bibliographystyle{elsarticle-num}
\bibliography{references,references_andrew,references_gregor,references_hendrik,references_trixi,references_local}

\appendix

\section{Algorithms and implementation}
\label{app:algorithms_implementation}

\linelabel{lne:appendix_start}
\addrev{3}{In the following, we present in detail the control flow for a single-physics simulation with
\trixi. The algorithms to enable coupled multi-physics simulations were already introduced in
\secref{multiphysics_coupling}. Here, the
purpose is to illustrate how simple it is to extend an existing single-physics solver for hyperbolic conservation
laws to a multi-physics solver for self-gravitating gas dynamics simulations.}

\addrev{3}{A two-dimensional quadtree mesh forms the basis of the simulation framework \trixi,
which can be refined adaptively during a simulation to meet dynamically changing resolution
requirements.  Several systems of equations are supported, including the compressible Euler
equations, ideal magnetohydrodynamics equations with divergence cleaning \cite{derigs2018ideal}, and
hyperbolic diffusion equations \cite{Nishikawa07,lou2019reconstructed}. They are discretized in
space by a high-order DG method and integrated in time by explicit Runge-Kutta schemes, as described
in \secref{model_methods}. \trixi is parallelized with a thread-based shared memory approach and
written in Julia \cite{bezanson2017julia}.}

\addrev{3}{During the \emph{initialization} phase, the simulation is set up
by creating a \emph{mesh} instance according to the simulation parameters.  Next, a \emph{solver}
instance is created that uses the mesh to build up the required data structures for the DG method,
and the solution is initialized. To obtain the optimal mesh for the initial solution state, the
adaptive mesh refinement algorithm (AMR) is called to adapt the mesh according to the initial
conditions, before initializing the solution again on this refined mesh.  This process is repeated
until the mesh no longer changes or a predetermined number of sub-cycles is reached.}

\addrev{3}{The execution phase begins when entering the \emph{main loop}, where first the current stable time
step is determined from the CFL condition. Next, the \emph{time integration} loop is entered, where
the solver is used to compute the spatial derivative once for each explicit Runge-Kutta stage. After
the solution has been advanced to the new time step, the AMR algorithm is called to adapt the mesh.
The main loop will then continue until either the desired simulation time is reached (flow-type
simulations) or the residual falls below a specified threshold that defines ``steady state''
(hyperbolic diffusion-type simulations).}

\addrev{3}{To adaptively refine the mesh, the adaptation algorithm first queries the DG solver for a refinement
indicator value $\lambda \in \{-1, 0, 1\}$ for each DG element.
The specific procedure to calculate $\lambda$ is problem dependent (one example is given in
\secref{sedov}).
An indicator value greater than zero
indicates that an element (and its corresponding grid cell) should be refined and a value less than
zero indicates that an element should be coarsened. A zero value means that the element should
remain unchanged. Next, the mesh and the solver are jointly adapted: In the first step, all cells
with an indicator value of $\lambda = 1$ are refined. A smoothing algorithm ensures that the
quadtree remains balanced, i.e., additional cells might be refined to retain a 2:1 relationship
between neighboring cells. The corresponding DG elements are then refined as well, using polynomial
interpolation to transfer the solution from the original coarse cell to the newly created refined
cell. In the second step, the same procedure is repeated to coarsen all cells with an indicator
value of $\lambda = -1$. In case of conflicts, refinement overrides coarsening to retain a high
solution quality. That is, if a cell originally marked for coarsening has already been refined in
the first step to rebalance the quadtree, it will not be coarsened.}
\linelabel{lne:appendix_end}


\end{document}